\documentclass[review,hidelinks,onefignum,onetabnum]{siamart250211}
\usepackage{enumitem}
\usepackage{multicol}

\usepackage{lipsum}
\usepackage{amsfonts,amssymb}
\usepackage{graphicx}
\usepackage{subfig}
\usepackage{epstopdf}
\usepackage{algorithmic}
\ifpdf
  \DeclareGraphicsExtensions{.eps,.pdf,.png,.jpg}
\else
  \DeclareGraphicsExtensions{.eps}
\fi

\newsiamremark{remark}{Remark}
\newsiamremark{hypothesis}{Hypothesis}
\crefname{hypothesis}{Hypothesis}{Hypotheses}
\newsiamthm{claim}{Claim}

\headers{A continuous approach to computing operator pseudospectra}{Kuan Deng, Xiaolin Liu, and Kuan Xu}

\title{A continuous approach to computing the pseudospectra of linear operators\thanks{Submitted to the editors DATE.
}}

\author{Kuan Deng\thanks{School of Mathematical Sciences, University of Science and Technology of China, 96 Jinzhai Road, Hefei 230026, Anhui, China (\email{dengkuan@mail.ustc.edu.cn}, \email{xiaolin@mail.ustc.edu.cn}, \email{kuanxu@ustc.edu.cn}).} \and Xiaolin Liu\footnotemark[2] \and Kuan Xu\footnotemark[2]}

\usepackage{amsopn}

\ifpdf
\hypersetup{
  pdftitle={A continuous approach to computing the pseudospectra of operators},
  pdfauthor={Kuan Deng, Xiaolin Liu, and Kuan Xu}
}
\fi

\def\mA{\mathcal{A}}
\def\mB{\mathcal{B}}
\def\mC{\mathcal{C}}
\def\mD{\mathcal{D}}
\def\mE{\mathcal{E}}
\def\mF{\mathcal{F}}
\def\mG{\mathcal{G}}
\def\mH{\mathcal{H}}
\def\mI{\mathcal{I}}
\def\mJ{\mathcal{J}}
\def\mK{\mathcal{K}}
\def\mL{\mathcal{L}}

\def\mS{\mathcal{S}}
\def\mP{\mathcal{P}}
\def\mQ{\mathcal{Q}}
\def\mO{\mathcal{O}}
\def\mR{\mathcal{R}}
\def\mT{\mathcal{T}}
\def\mU{\mathcal{U}}
\def\mV{\mathcal{V}}

\def\md{\mathrm{d}}

\def\bf{\boldsymbol{f}}

\def\hv{\hat{v}}

\def\mfr{\mathbf{r}}
\def\mfu{\mathbf{u}}
\def\mfv{\mathbf{v}}
\def\mfw{\mathbf{w}}
\def\mfx{\mathbf{x}}
\def\mfy{\mathbf{y}}

\def\tv{\tilde{v}}
\def\tw{\tilde{w}}

\def\tbv{\tilde{\mfv}}

\def\hbw{\hat{\mfw}}

\graphicspath{{./}}

\begin{document}

\maketitle


\begin{abstract}
We propose a continuous approach to computing the pseudospectra of linear operators with compact or compact-plus-scalar resolvent, following a ``solve-then-discretize'' strategy. Instead of taking a finite section approach or using a finite-dimensional matrix to approximate the operator of interest, the new method employs an operator analogue of the Lanczos process to work with operators and functions directly. The method is shown to be free of spectral pollution and spectral invisibility, fully adaptive, and nearly optimal in accuracy. The advantages of the method are demonstrated by extensive numerical examples and comparison with the traditional method.
\end{abstract}



\section{Introduction}\label{sec:intro}
In his 1999 Acta Numerica paper \cite{tre4}, Nick Trefethen wrote,

\bigskip

\textit{This field (the computation of pseudospectra) will also participate in a broader trend in the scientific computing of the future, the gradual breaking down of the walls between the two steps of discretization (operator $\longrightarrow$ matrix) and solution (matrix $\longrightarrow$ spectrum or pseudospectra).}

\bigskip

Not only serving as a projection for the future, this sentence also succinctly encapsulates the standard routine for computing the pseudospectra of an operator some 25 years ago. At the time of writing, it is still the default approach to the task. This approach begins by approximating a linear operator $\mL$, e.g., a differential operator, using a discretization method, e.g., a spectral method. If we denote by $L$ this finite-dimensional matrix approximation to $\mL$, the pseudospectra of $\mL$ at $z \in \mathbb{C}$ is then approximated by that of $L$ following the definition of the $\varepsilon$-pseudospectrum of a matrix
\begin{align}
\sigma_{\varepsilon}(L)=\left\{z \in \mathbb{C} \mid \Vert (zI - L)^{-1} \Vert >\varepsilon^{-1}  \right\} \label{mps},
\end{align}
where $I$ is the identity. In case of $z$ being an eigenvalue of $L$, we take the convention that $\Vert (z\mI-L)^{-1} \Vert = \infty$. If $2$-norm is adopted, \cref{mps} amounts to the computation of the smallest singular value of $zI - L$
\begin{align*}
\sigma_{\varepsilon}(L) = \{z \in \mathbb{C} \mid s_{\min}(zI - L) < \varepsilon\}.
\end{align*}
The standard method for calculating the pseudospectra of a matrix is the \textit{core EigTool algorithm} \cite[\S 39]{tre5}, which is recapitulated in \cref{alg:core}.

\begin{algorithm}
\renewcommand{\algorithmicrequire}{\textbf{Input:}}
\renewcommand{\algorithmicensure}{\textbf{Output:}}
\caption{Core EigTool algorithm}
\label{alg:core}
\begin{algorithmic}[1]
\REQUIRE{A matrix $L \in\mathbb{C}^{n\times n}$, $z\in \mathbb{C}$, a normalized vector $u_1 \in \mathbb{C}^n$, and a tolerance $\delta$}.
\STATE{Compute the Schur form $\tilde{L}$ so that $L$ and $\tilde{L}$ are unitarily similar.}
\REPEAT
\STATE{Perform Lanczos iteration $T(z)U_k = U_kH_k+\beta_{k+1}u_{k+1}e_k^T$, where $T(z)=(zI-\tilde{L})^{-*}(zI-\tilde{L})^{-1}$, $U_k=[u_1,\dots, u_k]\in\mathbb{C}^{n\times k}$, and $u_{k+1}\in \mathbb{C}^n$, to obtain the symmetric tridiagonal matrix $H_k\in \mathbb{R}^{k\times k}$.}
\STATE{Compute the largest Ritz value $\mu_1^{(k)} = \lambda_1(H_k)$.}
\STATE{$k \leftarrow k+1$}
\UNTIL{$| \mu_1^{(k)}/\mu_1^{(k-1)} - 1 | < \delta$}
\ENSURE{$\lVert (zI - L)^{-1} \rVert = \sqrt{\mu_1^{(k)}}$.}
\end{algorithmic}
\end{algorithm}


This ``discretize-then-solve'' paradigm has a few drawbacks: (1) There is no guarantee that the pseudospectra of $L$ serves as a quality approximation to that of the original operator, since $L$ is only of finite dimension. More importantly, the computation may suffer from the so-called \emph{spectral pollution}, i.e., spurious eigenvalues and pseudospectra, and \emph{spectral invisibility}, i.e., missing parts of the spectrum and pseudospectra. See \cref{sec:first} for an example. This is known as the ``finite section'' caveat \cite{han}. (2) It is difficult to effect adaptivity. (3) $T(z)$ might be more ill-conditioned than the original problem. See \cref{sec:advdiff} for an example. 
(4) Even though spectral approximation to the operators are used, the convergence of the resolvent norm may fall short of spectral in terms of the discretization size \cite[\S 43]{tre5}. See also \cref{sec:advdiff} for an example.

More recently, computing the pseudospectra of general bounded infinite matrices is addressed in \cite{han}, followed by a significantly improved version \cite{col3, col1} that offers a rigorous error control and the ability to deal with a greater range of problems, such as partial differential operators in unbounded domains. It is proposed to represent an operator $\mL$ as a matrix of infinite dimension that acts on $\ell^2(\mathbb{N})$. For a sequence of nested operator foldings and rectangular truncations of this infinite matrix, the limit of the smallest singular value is shown to be $\lVert(z\mI - \mL)^{-1}\rVert^{-1}$. This iterative method nevertheless also lacks adaptivity in stopping criterion and may have difficulty in dealing with operators restricted by boundary conditions, e.g., differential operators in compact domains. 

This paper shows that the pseudospectra of linear operators with compact or compact-plus-scalar resolvent can be computed following a ``solve-then-discretize'' strategy \cite{col4, col2, col5, gil1, hor1} by taking advantage of the recent development in spectral methods. The new method is an operator analogue of the Lanczos process that computes with operators and functions directly. This leads to a number of advantages: (1) We use effectively infinite-dimensional matrix representation of the operators, which rules out spectral pollution and spectral invisibility. (2) Justified by a priori analysis, both the application of resolvents and the stopping criterion of the Lanczos process are adaptive so that a nearly optimal accuracy can be obtained if a stringent exit condition is used. (3) We use the recently developed fast and well-conditioned spectral methods, making the most of the floating-point precision. (4) No involvement of quadrature or weight matrices, and the convergence is spectral. More importantly, the new method offers a unified framework for computing the pseudospectra of linear operators with compact or compact-plus-scalar resolvent. Though we confine ourselves in this article with certain spectral methods only, any method that solves the inverse resolvent equations of the operator of interest (see \cref{Eqn}) can be incorporated into the proposed framework straightforwardly in a plug-and-play fashion.

Throughout this paper, we assume $\mL: \mD(\mL) \subseteq \mH \rightarrow \mH $ be a closed linear operator on a separable Hilbert space $\mH$. For $z \in \mathbb{C}$, the resolvent of $\mL$ at $z$ is the operator
\begin{align*}
\mR(z) = (z\mI-\mL)^{-1},
\end{align*}
given that $\mR(z)\in \mB(\mH)$ exists, where $\mB(\mH)$ denotes the set of bounded operators on $\mH$. 
For $\varepsilon>0$, the $\varepsilon$-pseudospectrum of $\mL$ is
\begin{align}
\sigma_{\varepsilon}(\mL)=\left\{z \in \mathbb{C} \mid \Vert (z\mI-\mL)^{-1} \Vert >\varepsilon^{-1}  \right\} \label{epsps}
\end{align}
with the convention $\Vert (z\mI-\mL)^{-1} \Vert = \infty$ if $z \in \sigma(\mL)$, where $\sigma(\mL)$ denotes the spectrum of $\mL$. This definition can be found in, for example, \cite[\S 4]{tre5}. 

In the remainder of this paper, $\Vert \cdot \Vert = \Vert \cdot \Vert_2$, induced by the Euclidean inner product, unless indicated otherwise. We use $\langle \cdot, \cdot \rangle$ to denote the inner product. $\Lambda(\cdot)$ and $\sigma(\cdot)$ are used to denote the set of eigenvalues and the spectrum of an operator respectively. The largest eigenvalue, if exists and is isolated, is denoted by $\lambda_1(\cdot)$. $\Re(\cdot)$ stands for the real part of a complex number. $\mI$ and $I_{\infty}$ are respectively the identity operator and infinite-dimensional identity matrix. An asterisk is used to denote the complex conjugate of a complex variable or the adjoint of an operator. We denote by $\epsilon_{m}$ the machine epsilon, which is, for example, about $2.22 \times 10^{-16}$ in the double precision floating-point arithmetic. Throughout we follow \cite[p. 68]{hig} to take $\gamma_k = ck\epsilon_{m}/(1 - ck\epsilon_{m})$ for an integer $k$, where $c$ is a small integer constant whose value is unimportant. 

The paper is structured as follows. In \cref{sec:operator}, we introduce the operator analogue of the core EigTool algorithm. \Cref{sec:implementation} shows how the application of the resolvent is effected using Legendre-based spectral methods for adaptivity in resolution and dispensability of quadratures. In \cref{sec:conv&error}, we discuss the convergence and supply a careful error analysis which leads to an adaptive stopping criterion for the continuous Lanczos process. We test the method in \cref{sec:experiments} with extensive experiments before demonstrating a few important extensions in \cref{sec:extensions}. \Cref{sec:closing} closes by a summary.

\section{Operator analogue of the core EigTool algorithm}\label{sec:operator}

The definition \cref{epsps} offers limited guidance on the practical computation of an operator's pseudospectra. To relate the norm in \cref{epsps} to a computable quantity, we note that 
\begin{align*}
\lVert (z\mI-\mL)^{-1} \rVert^2 = \lVert \mR^*(z)\mR(z) \rVert = \sup\sigma(\mR^*(z)\mR(z)).
\end{align*}
If $\lVert \mR^*(z)\mR(z) \rVert$ happens to be the largest eigenvalue and is isolated, i.e, 
\begin{align}
\sup\sigma(\mR^*(z)\mR(z)) = \lambda_1(\mR^*(z)\mR(z)), \label{isolated}
\end{align}
\cref{alg:core} can be generalized to deal with $\mL$. \cref{alg:operator} is the operator analogue of the core EigTool algorithm, where the computation of the pseudospectra of an operator $\mL$ boils down to finding the largest eigenvalue of $\mT(z) = \mR^*(z)\mR(z)$.

\begin{algorithm}
   \renewcommand{\algorithmicrequire}{\textbf{Input:}}
   \renewcommand{\algorithmicensure}{\textbf{Output:}}
   \caption{The operator analogue of the core EigTool algorithm}\label{alg:operator}
   \begin{algorithmic}[1]
   \REQUIRE{A linear operator $\mL:\mD(\mL)\rightarrow \mH$, $z\in \mathbb{C}$, a normalized function $u_1\in \mH$, a tolerance $\delta$}
   \REPEAT
   \STATE{Perform $k$th Lanczos iteration (\cref{alg:lanczos})
   \begin{align}
   \mT(z)\mU_k=\mU_kH_k+\beta_{k+1}u_{k+1}e_k^T, \label{lanczos}
   \end{align}
   where $\mT(z)=(z\mI-\mL)^{-*}(z\mI-\mL)^{-1}$, $\mU_k = (u_1\vert \cdots \vert u_k)$, and $u_1, \dots, u_{k+1}\in \mH$, to obtain the symmetric tridiagonal matrix $H_k\in \mathbb{R}^{k\times k}$.}
   \STATE{Compute the largest Ritz value $\mu_1^{(k)} = \lambda_1(H_k)$ and the corresponding eigenvector $y_1^{(k)} = [y_{11}^{(k)},\dots,y_{k1}^{(k)}]^T\in \mathbb{R}^n$.}
   \STATE{$k \leftarrow k+1$}
   \UNTIL{$\beta_{k+1}|y_{k1}^{(k)}| <\max\left(C_L\epsilon_{m}\left(\mu_1^{(k)}\right)^{3/2},\ \delta\mu_1^{(k)}\right)$}
   \ENSURE{$\lVert (z\mI-\mL)^{-1} \rVert = \sqrt{\mu_1^{(k)}}$.}
   \end{algorithmic}
\end{algorithm}

The most notable distinction of \cref{alg:operator} from \cref{alg:core} lies in the substitution of matrix $L$ with the operator $\mL$. Additionally, there are two minor changes. First, the preliminary triangularization (line 1 of \cref{alg:core}), which is instrumental in reducing the computational cost \cite{lui}, has been removed. Second, a new stopping criterion is introduced in line 5 and further elaborated upon in \cref{sec:exit}.

Now we wonder for what kind of linear operator $\mL$ \cref{isolated} holds so that the pseudospectra of $\mL$ can be computed by \cref{alg:operator}. The following lemma shows that this is the case, if the resolvent $\mR(z)$ is compact or compact-plus-scalar.
\begin{lemma}\label{thm:compact}
(1) If $\mR(z)$ is compact, $\sup\sigma(\mT(z)) = \lambda_1(\mT(z))$. (2) If $\mR(z)$ is compact-plus-scalar, i.e., $\mR(z) = \mK + \tilde{z} \mI$, where $\tilde{z}\in \mathbb{C}$ and $\mK$ is compact on $\mH$, then $\sup\sigma(\mT(z)) = \lambda_1(\mT(z))$, except the case of $\sup\sigma(\mT(z)) = |\tilde{z}|^2$.
\end{lemma}

\begin{proof}
Statement (1) follows from the fact that $\mT(z) = \mR(z)^*\mR(z)$ is a compact self-adjoint  operator and the spectral theorem \cite[\S 4.5]{goh1}. 
To show (2), let
\begin{align*}
\mT(z) = (\mK + \tilde{z} \mI)^*(\mK + \tilde{z} \mI) = \mJ + |\tilde{z}|^2 \mI,
\end{align*}
where $\mJ = \tilde{z}\mK^* + \tilde{z}^*\mK + \mK^*\mK$. Noting that $\mJ$ is a compact self-adjoint operator, we use the spectral theorem again to have 
\begin{align*}
\sup\sigma(\mT(z)) = \sup\Lambda(\mJ) + |\tilde{z}|^2,
\end{align*}
where $\Lambda(\mJ)$ forms a real countably infinite sequence with a unique accumulation point at $0$. When $\sup\Lambda(\mJ) \neq 0$, $\sup\sigma(\mT(z)) = \lambda_1(\mJ) + |\tilde{z}|^2 = \lambda_1(\mT(z))$.
\end{proof}
The results we collect or show in the following subsections reveal that virtually all of the common linear operators have compact or compact-plus-scalar resolvent. These include differential operators restricted by proper boundary conditions and integral operators of Volterra and Fredholm types, which we now discuss individually. In addition, we also show how $\mR^*(z)$ is formulated for each operator.

\subsection{Differential operators}\label{sec:diff}
For the differential expression
\begin{align}
\tau = \frac{\md^N}{\md x^N} + a_{N-1}(x)\frac{\md^{N-1}}{\md x^{N-1}} + \cdots + a_1(x)\frac{\md}{\md x} + a_0(x), \label{diff}
\end{align}
where $a_k(x)$ for $k = 0, 1, \ldots, N-1$ are locally integrable on $[a, b]$, the \textit{maximal operator} $\mQ_{\max}: L^2([a,b]) \rightarrow L^2([a,b])$ is defined by $\mQ_{\max} g = \tau g$
with $\mD\left(\mQ_{\max}\right) = \{g \in \mA\mC_N([a,b]) \cap L^2([a,b]) \mid \tau g \in L^2([a,b])\}$\footnote{$\mA\mC_N([a,b])$ denotes the set of functions on $[a,b]$ whose $(N-1)$th derivative exists and is absolutely continuous.}. Suppose that $\mQ$ is a restriction of $\mQ_{\max}$ by proper boundary conditions \cite[\S XIV.3]{goh2}. It can be shown that $\mR(z) = (z\mI-\mQ)^{-1}$ is compact for $z \in \mathbb{C}\backslash \sigma(\mQ)$ \cite[\S XIV.3]{goh2}. For $\mR^*(z)$, we use the fact $(z\mI - \mQ)^{-*} = (z^*\mI - \mQ^*)^{-1}$, where $\mQ^*$ can be found by definition. See, for example, \cite[\S XIV.4]{goh2}.

\subsection{Fredholm integral operator}\label{sec:fred}
The Fredholm integral operator $\mF: L^2([a,b]) \rightarrow L^2([a,b])$ is defined as
\begin{align}
(\mF u)(s) = \int_{a}^{b} K(s,t)u(t) \md t,\quad a \le s \le b, \label{fred}
\end{align}
where the kernel $K(s,t)\in L^2([a,b]\times[a,b])$. For any nonzero $z\in \mathbb{C}$ that is not an eigenvalue of $\mF$, $\mR(z) = (z\mI-\mF)^{-1}$ is compact-plus-scalar \cite[Theorem 2.1.2]{smi1}. To evaluate $\mR^*(z)$, we use $(z\mI-\mF)^{-*} = (z^*\mI-\mF^*)^{-1}$, where $(\mF^*u)(s) = \int_a^b K^*(t,s)u(t)\md t$.


One of the most important subcases of the Fredholm integral operator is the Fredholm convolution integral operator
\begin{align}
(\mF_c u)(s) = \int_{a}^{b} K(s-t)u(t) \md t,\quad a \le s \le b. \label{fredconv}
\end{align}

\subsection{Volterra integral operator}\label{sec:volt}
The Volterra integral operator $\mV:L^2([a,b]) \rightarrow L^2([a,b])$ given by
\begin{align}
(\mV u)(s) = \int_{a}^{s}K(s,t)u(t)\md t,\quad a \le s \le b \label{volt}
\end{align}
can be deemed as a subcase of the Fredholm integral operator \cite[\S 2.7]{smi1} with the kernel $K(s,t)\in L^2([a,b]\times[a,b])$ zero-valued for $a \le s < t \le b$. For any nonzero $z \in \mathbb{C}$, $\mR(z) = (z\mI - \mV)^{-1}$ is compact-plus-scalar \cite[Theorem 2.7.1]{smi1}. For $\mR^*(z)$, we use the fact that $(z\mI-\mV)^{-*} = (z^*\mI-\mV^*)^{-1}$, where $(\mV^*u)(s) = \int_s^b K^*(t,s)u(t) \md t$. Similarly, an important subcase of the Volterra integral operator is the Volterra convolution integral operator
\begin{align}
(\mV_c u)(s) = \int_{a}^{s} K(s-t)u(t) \md t, \quad a \le t \le s \le b. \label{voltconv}
\end{align}



\subsection{Generalized eigenvalue problem}\label{sec:gep}
For the generalized eigenvalue problem (GEP)
\begin{align}
\mA x = \lambda \mB x, \label{gep}
\end{align}
we consider the case where both $\mA$ and $\mB$ are differential operators defined in \cref{sec:diff} subject to proper boundary conditions so that $\mB$ is invertible and $\mD(\mA) \subset \mD(\mB)$. Specifically, we assume that $\mA$ and $\mB$ have respectively the differential expressions
\begin{align}
\begin{aligned}
\tau_{\mA} &= \frac{\md^N}{\md x^N} + a_{N-1}(x)\frac{\md^{N-1}}{\md x^{N-1}} + \cdots + a_1(x)\frac{\md}{\md x} + a_0(x), \\
\tau_{\mB} &= b_M(x)\frac{\md^M}{\md x^M} + b_{M-1}(x)\frac{\md^{M-1}}{\md x^{M-1}} + \cdots + b_1(x)\frac{\md}{\md x} + b_0(x). \label{tauAB}
\end{aligned}
\end{align}
Without loss of generality, we also assume $M<N$. If we denote by $\Lambda(\mA, \mB) = \{\lambda\in \mathbb{C} \mid \mA x = \lambda \mB x,\ x\ne 0 \}$ the set of eigenvalues of \cref{gep}, \cref{gep} can be solved as the standard eigenvalue problem with $\Lambda(\mA, \mB) = \Lambda(\mG)$ by the invertibility of $\mB$, where $\mG = \mB^{-1}\mA$. Moreover, this is followed by the $\varepsilon$-pseudospectrum of the GEP
\begin{align*}
\sigma_{\varepsilon}(\mA, \mB)=\sigma_{\varepsilon}\left(\mG \right),
\end{align*}
which is the operator extension of the definition for the pseudospectra of matrix GEPs \cite[\S 45]{tre5}.  


\begin{lemma}\label{lem:Gep} For the GEP \cref{gep}, we further assume that $a_j(x), b_k(x) \in L^{\infty}([a,b])$ for $j = 0, 1, \ldots, N-1$ and $k = 0, 1, \ldots, M$ respectively and $\mD(\mA^*) \subset \mD(\mB^*)$. Then $\mR(z) = (z\mI-\mG)^{-1}$ and 
\begin{subequations}
\begin{align}
\mR^*(z) = \mB^*(z^*\mB^*-\mA^*)^{-1} \label{R*}
\end{align}
are both compact for $z \in \mathbb{C}\backslash \sigma(\mG)$; 
In addition, the restriction of $\mR(z)$ to $\mD(\mB_{\max})$ is 
\begin{align}
\mR(z)\big|_{\mD(\mB_{\max})} = (z\mB-\mA)^{-1}\mB_{\max}. \label{R} 
\end{align}
\end{subequations}
\end{lemma}
\begin{proof}
We note that $(z\mI-\mG)^* = (z^*\mB^*-\mA^*)\mB^{-*}$, and its domain $\mD((z\mI-\mG)^*) = \{u \in L^2([a,b])\arrowvert\mB^{-*}u \in \mD(\mA^*)\}$. When $z^*\mB^*-\mA^*$ is invertible, so is $(z\mI-\mG)^*$,
and $\mR^*(z) =  \left((z\mI-\mG)^{-1}\right)^* = \left((z\mI-\mG)^*\right)^{-1} = \mB^*(z^*\mB^{*}-\mA^{*})^{-1}$ by Proposition 2.6 in \cite[\S XIV.2]{goh2}. Since $\{a_k\}_{k = 0}^{N-1}, \{b_k\}_{k = 0}^M \in L^{\infty}([a,b])$ and $M<N$, $\mB^*$ is $(z^*\mB^{*}-\mA^{*})$-compact by Theorem 1.3 in \cite[\S XVII.1]{goh2}. It then follows from the invertibility of $z^*\mB^{*}-\mA^{*}$ and the argument used in the proof of Theorem 5.1 in \cite[\S XVII.5]{goh2} that $\mB^*(z^*\mB^{*}-\mA^{*})^{-1}$ is compact, i.e., $\mR^*(z)$ is compact, which, in turn, leads to the compactness of $\mR(z)$.

To show \cref{R} , we first consider $\mR(z)\big|_{\mD(\mB)}$. For $\phi \in \mD(\mB)$ and $\psi \in L^2([a,b])$, 
\begin{align}
\langle (z\mB-\mA)^{-1}\mB \phi, \psi\rangle = \langle \mB \phi, (z^*\mB^{*}-\mA^{*})^{-1}\psi\rangle = \langle \phi, \mB^{*}(z^*\mB^{*}-\mA^{*})^{-1}\psi\rangle, \label{inner}
\end{align}
from which it follows that $\mR(z)\big|_{\mD(\mB)} = (z\mB-\mA)^{-1}\mB$. By the extension theorem for bounded linear operators, we have $\mR(z)\big|_{\mD(\mB_{\max})} = (z\mB-\mA)^{-1}\mB_{\max}$.
\end{proof}

When evaluate $\mR(z)$ and $\mR^*(z)$ by \cref{R} and \cref{R*} respectively, we need to know $\mA^*$ and $\mB^*$, which can be figured out as for the differential operators in \cref{sec:diff}.

\section{Implementation}\label{sec:implementation}
In this section, we focus on the practical implementation of lines 2, i.e., the Lanczos iteration. We defer the discussion on the stopping criterion to \cref{sec:exit}. To facilitate the discussion, we elaborate the Lanczos iteration in \cref{alg:lanczos}, which can be embedded into \cref{alg:operator} at line 2 there. 

\begin{algorithm}
\begin{multicols}{2}
\renewcommand{\algorithmicrequire}{\textbf{Input:}}
\renewcommand{\algorithmicensure}{\textbf{Output:}}
\caption{Lanczos process}\label{alg:lanczos}
\begin{algorithmic}[1]
\REQUIRE{$\mT(z)$, $H_{k-1}$, $\beta_{k}$, $u_{k}$, $u_{k-1}$}
\STATE{$w_k = \mT(z)u_k$}
\IF{$k>1$}
\STATE{$p_k = w_k - \beta_k u_{k-1}$}
\ELSE 
\STATE{$p_k = w_k$}
\ENDIF
\STATE{$\alpha_k = \Re(\langle u_k, p_k \rangle)$}
\STATE{$q_k = p_k - \alpha_k u_k$}
\STATE{$\beta_{k+1} = \lVert q_k \rVert$}
\STATE{$u_{k+1} = q_k/\beta_{k+1}$}
\ENSURE{$H_k$, $\beta_{k+1}$, $u_{k+1}$}
\end{algorithmic}
\end{multicols}
\end{algorithm}

In line 1 of \cref{alg:lanczos}, we need to apply $\mT(z)=(z\mI-\mL)^{-*}(z\mI-\mL)^{-1}$ to $u \in \mH$. This amounts to solving 
\begin{subequations}
\begin{align}
(z\mI-\mL)v &= u,  \label{resEqn} \\
(z^*\mI-\mL^*) w &= v \label{adjEqn}
\end{align}\label{Eqn}%
\end{subequations}
in sequence. We solve these operator equations using coefficient-based spectral methods with normalized Legendre polynomials. For instance, for \cref{resEqn} the solution $v$ and the right-hand side $u$ are represented by series expansion in normalized Legendre polynomials. This choice is motivated by the fact that normalized Legendre polynomials form an orthonormal basis of $L^2([-1,1])$. As a result, the $L^2$-norm of a series of normalized Legendre polynomials coincides exactly with the $l^2$-norm of the series coefficients, and the $L^2$ inner product of two series of normalized Legendre polynomials corresponds to the dot product of the coefficient vectors. These identities eliminate the need for polynomial multiplication and numerical quadrature when computing inner products and norms in the Lanczos process; see lines 7 and 9 in \cref{alg:lanczos}.


The variable coefficients in a differential operator and the kernel function of an integral operator are assumed to have certain regularity so that they can be approximated by a finite normalized Legendre series or a low-rank adaptive cross approximant \cite{tow1} up to any desired accuracy. The coefficients of these series and adaptive cross approximants can be obtained using the fast Legendre transform \cite{kei} with the degrees of the series determined by, e.g., the chopping algorithm \cite{aur1}. In this section and \cref{sec:conv&error}, an operator carrying a tilde is an approximation to the original one obtained by replacing the variable coefficients or kernel function by finite-degree polynomials, and any tilde-decorated variable is the perturbed version due to the adoption of the approximate operator. 

To make the presentation self-contained, we now give the detail of the spectral methods for each operator by first focusing on constructing the discrete linear system for \cref{Eqn} for each operator discussed in \cref{sec:operator}. Except for the GEP, we discuss only \cref{resEqn}, and analogue is easily drawn for \cref{adjEqn}. 
Since all these spectral methods lead to infinite-dimensional banded systems, we defer the discussion of the solution technique to \cref{sec:solving}.


\subsection{Differential operator}\label{sec:diffimp} When $\mL$ is a differential operator $\mQ$ defined in \cref{sec:diff}, we solve \cref{Eqn} by the ultraspherical spectral method\footnote{The method used here is a variant of the original ultraspherical spectral method in that the latter works with Chebyshev polynomials and ultraspherical polynomials of integer orders, i.e., $C^{(N)}(x)$ for $N = 1, 2, \ldots$, whereas we employ the ultraspherical polynomials of half orders, i.e., $C^{(N+1/2)}(x)$ for $N = 0, 1, \ldots$, since the solution is represented in normalized Legendre polynomials.} \cite{olv1} with basis recombination \cite{qin}. Specifically, we approximate the variable coefficients $a_{N-1}(x),\ldots, a_0(x)$ by normalized Legendre series of degrees, say, $m^{a_{N-1}},\ldots, m^{a_0}$. Denoting by $\tilde{\mQ}$ the approximate differential operator, we aim at solving the approximate problem
\begin{align}
(z\mI - \tilde{\mQ}) \tv = u. \label{resEqn1}
\end{align}
Instead of representing the solution directly using the standard normalized Legendre polynomials, we first express the solution in a recombined basis formed by suitable combinations of neighboring normalized Legendre polynomials that satisfy the boundary conditions.\footnote{Such a basis is commonly referred to as a recombined basis, and is also known by various other names in different contexts, including integrated Legendre polynomials (of Babu\v{s}ka), bubble functions, modal basis, and boundary-adapted basis. The use of such recombined bases in spectral methods was pioneered by Shen \cite{she,she2}.} With this recombined basis, the matrix representation of $z\mI - \tilde{\mQ}$ becomes an infinite-dimensional strictly banded matrix, which can be efficiently factorized by subroutines for banded matrices in a standard library.

More importantly, the range of this matrix representation lies entirely in the infinite-dimensional vector space associated with the coefficients of the $C^{(N+1/2)}(x)$ expansion.\footnote{If the standard ultraspherical spectral method is applied to solve \cref{Eqn} without basis recombination, the range of the left-hand-side matrix would instead be the union of this space and a finite-dimensional space corresponding to the boundary conditions.} This property is crucial for establishing \cref{thm:v-vt}.

On the right-hand side, the normalized Legendre coefficients of $u$ must likewise be promoted to the ultraspherical basis $C^{(N+1/2)}(x)$ in order to match the left-hand side. After solving the resulting system, the solution, as the series coefficients in recombined normalized Legendre polynomials, is transformed back to those in the standard normalized Legendre polynomials.

\subsection{Fredholm operator}\label{sec:fredimp} When $\mL$ is a general Fredholm operator $\mF$ defined by \cref{fred}, we first approximate the kernel $K(s,t)$ by a low-rank adaptive cross approximation of rank $r$, i.e.,
\begin{align}
K(s,t) \approx \tilde{K}(s,t) = \sum_{j=1}^r f_j(s) g_j(t). \label{approx2}
\end{align}
Here, $f_j$ and $g_j$ $(j = 1, 2, \ldots, r)$ are the normalized Legendre series of degree $m^f$ and $m^g$ respectively. Thus, the approximate operator equation we actually solve is
\begin{align}
(z\mI - \tilde{\mF}) \tv = u, \label{resEqn2}
\end{align}
where $\tilde{\mF}$ is the Fredholm operator defined by the approximate kernel $\tilde{K}(s,t)$. 
Since
\begin{align*}
( z\mI - \tilde{\mF} ) \tv \approx z \tv(s) - \sum_{j=1}^r \langle g_j^*(t), \tv(t) \rangle f_j(s), \label{zIFv}
\end{align*}
the matrix representation of $z\mI - \tilde{\mF}$ is the sum of an infinite-dimensional diagonal and a rank-$k$ semiseparable matrix
\begin{align*}
M_{\mF}(z) = zI_{\infty} - CR^*,
\end{align*}
where $C = \left[C_1, \cdots, C_r \right] \in \mathbb{C}^{\infty \times r}$ and $R = \left[ R_1, \cdots, R_r \right] \in \mathbb{C}^{\infty \times r}$. The top $m^f+1$ elements of $C_j$ are the coefficients of normalized Legendre expansion of $f_j(s)$ and the rest of $C_j$ are zeros. Similarly, $R_j$ is zero except the first $m^g+1$ elements which store the normalized Legendre coefficients of $g_j(t)$. Thus, only the top left $(m^f+1) \times (m^g+1)$ block and the diagonal of $M_{\mF}$ are nonzero. 


\subsection{Fredholm convolution operator}\label{sec:freconvimp} When $\mL$ is a Fredholm convolution integral operator $\mF_c$ given by \cref{fredconv}, we approximate the univariate kernel function $K(x)$ by a normalized Legendre series $\tilde{K}(x)$ of degree, say, $m^K$, yielding an approximate convolution operator $\tilde{\mF}_c$. The approximate operator equation that we solve is
\begin{align}
(z\mI - \tilde{\mF_c}) \tv = u. \label{resEqn3}
\end{align}
Then the matrix representation of $\tilde{\mF}_c$ is an infinite-dimensional matrix with only the upper skew-triangular part of the top left $(m^K+1) \times (m^K+1)$ block being nonzero. These nonzero entries can be obtained from the $m^K+1$ coefficients of $\tilde{K}(x)$ using a four-term recurrence formula \cite{liu}. Adding it to the infinite-dimensional scalar matrix representing $z\mI$ gives the infinite-dimension arrow-shaped representation of $z\mI-\tilde{\mF}_c$.

\subsection{Volterra operator}\label{sec:voltimp} In case of $\mL$ being the Volterra operator $\mV$ defined by \cref{volt}, we follow the same strategy for handling the Fredholm operator by approximating the kernel $K(s,t)$ with the low-rank approximation \cref{approx2} so that we solve the resulting approximate operator equation
\begin{align}
(z\mI - \tilde{\mV}) \tv = u, \label{resEqn4}
\end{align}
where
\begin{align*}
(\tilde{\mV} \tv)(s) = \sum_{j=1}^{r} f_j(s)\int_a^s g_j(t) \tv(t)dt.
\end{align*}
The matrix representation of $z\mI-\tilde{\mV}$ is again an infinite-dimensional banded matrix
\begin{align*}
M_{\mV}(z) = zI_{\infty} - \sum_{j=1}^{r} M_{\frac{1}{2}}[f_j]J_{[a,b]} M_{\frac{1}{2}}[g_j], 
\end{align*}
where $J_{[a,b]}$ is the indefinite integral matrix for normalized Legendre series \cite[\S 4.7]{sze} and $M_{\frac{1}{2}}[f_j]$ and $M_{\frac{1}{2}}[g_j]$ are the normalized Legendre multiplication matrices \cite[\S 3.1]{olv1}.


\subsection{Volterra convolution operator}\label{sec:voltconvimp} When $\mL$ is a Volterra convolution integral operator $\mV_c$ defined by \cref{voltconv}, we approximate the kernel $K(x)$ by a normalized Legendre series of degree $m^K$. This gives the approximate operator equation
\begin{align}
(z\mI - \tilde{\mV}_c) \tv = u. \label{resEqn5}
\end{align}
With the coefficients of this finite normalized Legendre series, a banded infinite-dimensional matrix representation of $\tilde{\mV}_c$ can be constructed using a three-term recurrence formula and a scaled symmetry relation \cite{hal,xu}. Adding $zI_{\infty}$ does not change the sparsity structure.

\subsection{Generalized eigenvalue problem}\label{sec:gepimp} For the generalized eigenvalue problem defined by \cref{gep},
it follows from \cref{R} that $v = (z\mB-\mA)^{-1}\mB_{\max}u$, or equivalently 
\begin{align}
(z\mB-\mA) v = \mB_{\max}u, \label{Ru}
\end{align}
in place of \cref{resEqn}. To solve \cref{Ru} with the ultraspherical spectral method, we follow the routine to replace the variable coefficients in \cref{tauAB} by normalized Legendre series, as for the differential operators. This gives the approximate operator equation
\begin{align}
(z\tilde{\mB}-\tilde{\mA}) \tv = \tilde{\mB}_{\max}u, \label{resEqn6}
\end{align}
where $\tilde{\mA}$, $\tilde{\mB}$, and $\tilde{\mB}_{\max}$ are the approximations to $\mA$, $\mB$, and $\mB_{\max}$ respectively. Again, the matrix representation of the differential operator $z\tilde{\mB}-\tilde{\mA}$ can be made strictly banded with recombined normalized Legendre polynomials. On the right-hand side, $\tilde{\mB}_{\max}$ is also banded and its product with $u$ is promoted to the $C^{(N+1/2)}(x)$ space to match the left-hand side. 

Suppose that $\tv$ is available from solving \cref{resEqn6}. To proceed to the application of $\mR^*(z)$ to $\tv$ following \cref{R*}, 
we solve the approximate equation
\begin{align*}
(z^*\tilde{\mB}^*-\tilde{\mA}^*)\tilde{h} = \tv, \label{R*vt}
\end{align*}
which is obtained again by replacing the variable coefficients in 
$(z^*\mB^*-\mA^*)h = v$ with finite normalized Legendre series. Once $\tilde{h}$ is available, we calculate $\tw = \tilde{\mB}^*_{\max}\tilde{h}$ to recover the approximate version of $w$. This again amounts to the product of a banded matrix and a vector. Note that the notations $\tilde{\mA}^*$, $\tilde{\mB}^*$, and $\tilde{\mB}^*_{\max}$ are adopted for simplicity and readability; they do not necessarily equal the adjoints of $\tilde{\mA}$, $\tilde{\mB}$, and $\tilde{\mB}_{\max}$, respectively. 


\subsection{Solving the banded system}\label{sec:solving}
Since \cref{resEqn1,resEqn4,resEqn5,resEqn6} all lead to infinite-dimensional banded systems and \cref{resEqn2,resEqn3} give infinite-dimensional system that can be regarded as banded, we discuss their solution by looking at a unified model equation
\begin{align}
(z\mI - \tilde{\mL}) \tv = u. \label{resEqnt}
\end{align}
Here, $\tilde{\mL}$ is obtained, as detailed above, by replacing the exact variable coefficients or kernel function by the finite normalized Legendre series. Suppose that $\tilde{\mL}$ deviates from $\mL$ by $\Delta \mL$, i.e., $\mL = \tilde{\mL} + \Delta \mL$, and \cref{resEqnt} leads to the infinite-dimensional linear system  
\begin{align}
A \tbv = \mfu \label{system}
\end{align}
with the lower bandwidth of $A$ being $b$ and the first $n_u$ elements of $\mfu$ being nonzero. The infinite vectors $\tbv$ and $\mfu$ collect the coefficients of $\tv$ and $u$ in the bases $\{\phi_i\}_{i = 1}^{\infty}$ and $\{\psi_i\}_{i = 1}^{\infty}$ respectively. For differential operators or generalized eigenvalue problems, $\{\phi_i\}_{i = 1}^{\infty}$ and $\{\psi_i\}_{i = 1}^{\infty}$ are recombined normalized Legendre and ultraspherical polynomials respectively; for integral operators they are simply the standard normalized Legendre polynomials.

Following \cite{olv1}, we use the adaptive QR to triangularize $A$ with a sequence of Householder reflectors and applies the same reflectors to the right-hand side $\mfu$ simultaneously. For $n \geq n_u$, we denote by $Q_n\in \mathbb{C}^{(n+b) \times (n+b)}$ the product of the first $n$ Householder reflectors. After the first $n$ steps, the left-hand side of the partially triangularized system is
\begin{align*}
\begin{pmatrix}
Q_n^* & \phantom{0}\\
\phantom{0}  & I_{\infty}
\end{pmatrix}
\begin{pmatrix}
A_{11} & A_{12}\\
\phantom{0} & A_{22}
\end{pmatrix}
= 
\begin{pmatrix}
R_{11} & Q_n^* A_{12}\\
\phantom{0} & A_{22}
\end{pmatrix}, \label{adaQRA}
\end{align*}
where $A_{11} \in \mathbb{C}^{(n+b) \times n}$ and $R_{11} \in \mathbb{C}^{n \times n}$ is upper triangular. The resulting infinite vector $\mfr$ on the right-hand side is
\begin{align*}
\mfr =
\begin{pmatrix}
Q_n^* & \phantom{0}\\
\phantom{0}  & I_{\infty}
\end{pmatrix}
\begin{pmatrix}
\mfu_{1:n_u}\\
0
\end{pmatrix}
= 
\begin{pmatrix}
\mfr_{1:n} \\
\mfr_{n+1:n+b}\\
0
\end{pmatrix}, \label{adaQRu}
\end{align*}
whose first $n+b$ elements are nonzero. We can then obtain a truncated solution
\begin{align}
\bar{\mfv} = 
\begin{pmatrix}
R_{11}^{-1}\mfr_{1:n}\\
0
\end{pmatrix} \label{vbar}
\end{align}
by back substitution. In exact arithmetic, $\mfr_{n+1:n+b}$ is exactly the error of $\bar{\mfv}$ in residual \cite{olv1}. Thus, we terminate the QR process if $\lVert \mfr_{n+1:n+b} \rVert < \epsilon \lVert \mfr \rVert$ for a preset tolerance $\epsilon$. Otherwise, we keep applying Householder reflectors until $\lVert \mfr_{n+1:n+b} \rVert < \epsilon \lVert \mfr \rVert$ is satisfied. In the rest of this paper, we refer to this value of $n$ as the degrees of freedom (DoF). In the current investigation, we choose $\epsilon = \epsilon_{m}$ since we intend to make the most of the digits in floating-point arithmetic and achieve the highest possible accuracy. 


\section{Convergence, error analysis, and stopping criterion}\label{sec:conv&error}
In this section, we discuss the convergence of the proposed method (\cref{sec:conv}) briefly, followed by a careful error analysis (\cref{sec:error}) which leads to an adaptive stopping criterion and an estimate for the relative error in the computed resolvent norm (\cref{sec:exit}).

\subsection{Convergence} \label{sec:conv}
The convergence theory of the operator Lanczos iteration is a well known result due to Saad \cite{saad1}. The theorem that follows extends the original result of Saad for a compact self-adjoint operator to a self-adjoint operator that satisfies \cref{isolated}.

\begin{theorem}\label{thm:Convergence}(Convergence of Lanczos iteration) Consider a self-adjoint operator $\mS$ on a Hilbert space $\mH$. Suppose $\lambda_1(\mS) = \lVert \mS\rVert$ is the largest eigenvalue, corresponding to eigenmode $\varphi_1$.
Let $\lambda_{\inf} = \inf\sigma(\mS)$ represent the infimum of the spectrum of $\mS$. For $u_1$ that satisfies $P_1 u_1 \neq 0$, after $k$ Lanczos iterations
\begin{align}
0 \le \lambda_1 - \mu_1^{(k)} \le (\lambda_1-\lambda_{\inf})\frac{\tan ^2\theta(\varphi_1, u_1)}{T^2_{k-1}(\omega)}, \label{convergence}
\end{align}
where $\theta(\varphi_1, u_1)$ is the angle between $\varphi_1$ and $u_1$, $P_1$ is the orthogonal projection on the eigenspace corresponding to $\lambda_1$, $T_m$ is $m$th Chebyshev polynomial of the first kind\footnote{Since $\omega$ is greater than $1$, the Chebyshev polynomial $T_{k-1}$ in \cref{convergence} is evaluated outside the canonical domain $[-1, 1]$.}, and 
\begin{align*}
\omega = 1 + 2\frac{\lambda_1-\lambda_2}{\lambda_2-\lambda_{\inf}}, \quad \lambda_2 = \sup(\sigma(\mS)\backslash \lambda_1).  
\end{align*}
\end{theorem}

A little algebraic work shows that \cref{convergence} can be simplified as
\begin{align*}
|\lambda_1 - \mu_1^{(k)}| \le c \omega^{-k}, \label{exponential}
\end{align*}
where $c$ is a constant independent of $\omega$. This confirms what we observed---the Ritz values $\mu_1^{(k)}$ obtained in \cref{alg:operator} converges to $\lambda_1(\mT(z))$ exponentially fast. While \cref{thm:Convergence} holds in exact arithmetic, the effect of rounding error is nevertheless not negligible in floating-point arithmetic and is examined by the analysis below.

\subsection{Error analysis}\label{sec:error}
In this subsection, we derive an error bound for the computed resolvent norm obtained by \cref{alg:operator,alg:lanczos}. Our analysis consider both the mathematical errors inherent in the algorithms, e.g., due to truncation, and the numerical errors caused by rounding in floating-point arithmetic. In the remainder of this subsection, a hat over a variable indicates that this is the computed version in floating-point arithmetic of the unhatted variable.

We commence by revisiting the solution of \cref{system} via the adaptive QR method. There is of course no way to solve \cref{system} for the exact partial solution $\bar{\mfv}$ given in \cref{vbar} in floating-point arithmetic. Let $\hat{\mfv}$ and $\hat{\mfr}$ be the computed solution and the numerical version of $\mfr$ respectively after $n$ steps of QR factorization via Householder reflection. The lemma that follows gives the backward errors in the computed solution to \cref{system}, which will be useful for \cref{thm:v-vt}. Note that these errors incorporate both the mathematical and numerical ones. We denote by $\mP_n = (I_n, 0)$ the $n \times \infty$ projection matrix.

\begin{lemma}\label{lem:lempert}
Suppose that we solve \cref{system} by means of Householder reflection in floating-point arithmetic and stop after $n$ steps so that $\lVert \hat{\mfr}_{n+1:n+b} \rVert < \epsilon_{m} \lVert \hat{\mfr} \rVert$. It then follows
\begin{align}
(A+E)\hat{\mfv} &= \mfu - \Delta \mfu, \label{pertinfeq}  
\end{align}
where $E = \mP_{n+b}^{\top}E_{11}\mP_n$ and $E_{11} \in \mathbb{C}^{(n+b)\times n}$. In addition, the first $n+b$ elements of $\Delta \mfu$ are nonzero. These perturbations are bounded as
\begin{subequations}
\begin{align}
\lVert E \rVert &\leq \gamma_{bn^{3/2}} \lVert A_{11} \rVert, \label{normE} \\
\lVert \Delta \mfu\rVert &\leq \gamma_{bn} \lVert \mfu \rVert. \label{normdu}
\end{align}\label{pertnorm}%
\end{subequations}
\end{lemma}
\begin{proof}
It follows from \cite[Theorem 19.5]{hig} that
\begin{align*}
(A_{11} + E_{11}) \hat{\mfv}_{1:n} = Q_n
\begin{pmatrix}
\hat{\mfr}_{1:n}\\
0
\end{pmatrix} = Q_n
\begin{pmatrix}
\mfr_{1:n} + \Delta \mfr_{1:n} \\
0
\end{pmatrix},
\end{align*}
where $\lVert E_{11} \rVert \leq \gamma_{bn^{3/2}} \lVert A_{11} \rVert$ and $\lVert \Delta \mfr_{1:n} \rVert \leq \lVert \mfr - \hat{\mfr} \rVert \leq \gamma_{bn} \lVert \mfr \rVert$. The former is equivalent to \cref{normE}. Embedded into the infinite-dimensional linear system, the last equation amounts to 
\begin{align*}
(A + E) \hat{\mfv} = 
\begin{pmatrix}
Q_n& \phantom{0}\\
\phantom{0}  & I_{\infty}
\end{pmatrix} 
\begin{pmatrix}
\mfr_{1:n} + \Delta \mfr_{1:n} \\
0\\
0
\end{pmatrix},
\end{align*}
which is exactly \cref{pertinfeq} with
\begin{align*}
   \Delta \mfu = 
   \begin{pmatrix}
   Q_n& \phantom{0}\\
   \phantom{0}  & I_{\infty}
   \end{pmatrix} 
   \begin{pmatrix}
   -\Delta \mfr_{1:n}\\
   \hat{\mfr}_{n+1:n+b}\\
   0
   \end{pmatrix}.
\end{align*}
Since $Q_n$ is orthogonal and $\lVert \hat{\mfr}_{n+1:n+b} \rVert < \epsilon_{m} \lVert \hat{\mfr} \rVert$, $\lVert \Delta \mfu\rVert \leq \lVert \Delta \mfr_{1:n} \rVert + \lVert \hat{\mfr}_{n+1:n+b} \rVert= \gamma_{bn} \lVert \mfu \rVert$. 
\end{proof}


For $E$ in \cref{pertinfeq}, we now define
\begin{align*}
\mE : \mD(\tilde{\mL}) \rightarrow \mH,\quad \mE x = y, \quad \text{s.t.} \quad \mfy = E \mfx, \label{pertop}  
\end{align*}
where $\mfx$ and $\mfy$ are the coefficients of $x$ and $y$ in the bases $\{\phi_i\}_{i = 1}^{\infty}$ and $\{\psi_i\}_{i = 1}^{\infty}$ respectively. Thus, the finite-rank operator $\mE$ is the perturbation of $z\mI-\tilde{\mL}$. In order to bound $\mE$ by making use of \cref{normE}, we need to relate the norm of a polynomial in $\mD(\tilde{\mL})$ or $\mH$ to that of its coefficients in $\{\phi_i\}_{i = 1}^{\infty}$ or $\{\psi_i\}_{i = 1}^{\infty}$ respectively. In fact, for $x \in \Phi_n := span \{\phi_i\}_{i = 1}^n$ and $y \in  \Psi_n := span \{\psi_i\}_{i = 1}^n$ the norms of the vectors $\mfx$ and $\mfy$ are sandwiched by those of $x$ and $y$ as
\begin{subequations}
\begin{align}
&\alpha^{\phi}_n \lVert x \rVert \leq \lVert \mfx \rVert \leq \beta^{\phi}_n\lVert x \rVert, \quad \forall x \in \Phi_n, \label{basisL} \\
&\alpha^{\psi}_n \lVert y \rVert \leq \lVert \mfy \rVert \leq \beta^{\psi}_n\lVert y \rVert, \quad \forall y \in \Psi_n. \label{basisR}
\end{align}\label{basisCond}%
\end{subequations}
Here, the parameters $\alpha^{\phi}_n$, $\beta^{\phi}_n$, $\alpha^{\psi}_n$, and $\beta^{\psi}_n$ are basis-dependent. We also define the condition number of these bases as
\begin{align*}
\kappa^{\phi}_n = \frac{\beta^{\phi}_n}{\alpha^{\phi}_n}, \quad \kappa^{\psi}_n = \frac{\beta^{\psi}_n}{\alpha^{\psi}_n}.
\end{align*}

For differential operators, since $\phi_i$ is a combination of normalized Legendre polynomials that satisfies boundary conditions, the condition number $\kappa^{\phi}_n$ also depends on the boundary conditions. For example, $\kappa^{\phi}_n = O(n)$ for Dirichlet boundary condition and $\kappa^{\phi}_n = O(n^2)$ for clamped boundary conditions \cite[\S 4.1.5]{mor}. If the differential operator of interest is an $N$th-order one, $\{\psi_i\}_{i=0}^{\infty}$ are the ultraspherical polynomials $C^{(N+1/2)}(x)$ and the corresponding condition number $\kappa^{\psi}_n$ is then $O(n^N)$. Of course, these estimates are rather pessimistic in that the effective condition numbers in practice are much smaller than these asymptotic results. For integral operators, $\kappa^{\phi}_n=\kappa^{\psi}_n=1$ because $\{\phi_i\}_{i = 1}^{\infty}$ and $\{\psi_i\}_{i = 1}^{\infty}$ are simply the normalized Legendre polynomials.

Let $\hat{v}$ be the polynomial corresponding to $\hat{\mfv}$. With mild assumptions, the lemma below, where $\Delta \mL$, $E$, and $\Delta \mfu$ are all factored in, bounds the forward error in the computed solutions to \cref{resEqn}, i.e., $v-\hv$, in terms of $\xi = \lVert\mT(z)\rVert$. The definition \cref{epsps} suggests that $\xi$ is $\mO(\varepsilon^{-2})$. Since we are interested in $\varepsilon \ll 1$, $\xi$ is usually very large.


\begin{lemma}\label{thm:v-vt}
If $\lVert \mR(z) \Delta \mL \rVert < \lVert \mR(z) \rVert\gamma_1  < 1/2$ and $\lVert \mR(z) \rVert\lVert \mE \rVert < 1/2$, the forward error
\begin{align}
\lVert v - \hat{v} \rVert \leq (C_1\xi + C_2\xi^{1/2})\lVert u \rVert, \label{errorv}
\end{align}
where
\begin{align*} 
C_1 = C_3 \kappa^{\phi}_n & \kappa^{\psi}_{n+b} \lVert (z\mI -\mL)|_{\Phi_n}\rVert \gamma_{bn^{3/2}}, \quad C_2 = C_3\frac{\beta^{\psi}_{n_u}}{\alpha^{\psi}_{n+b}} \gamma_{bn},\\
 &C_3 = \frac{1}{1 - \lVert \mR(z) \Delta \mL \rVert - \lVert \mR(z) \rVert \lVert \mE \rVert }.
\end{align*}
\end{lemma}

\begin{proof}
We first note from \cref{pertinfeq} that
\begin{align}
(z\mI - \tilde{\mL} + \mE)\hat{v} = u - \Delta u, \label{opperteq}
\end{align}
where $\Delta u$ is the polynomial expressed in the basis $\{\psi_i\}_{i = 1}^{\infty}$ corresponding to $\Delta \mfu$. It follows from \cref{basisR} and \cref{normdu} that
\begin{align}
\lVert \Delta u \rVert \leq \frac{\beta^{\psi}_{n_u}}{\alpha^{\psi}_{n+b}}\gamma_{bn}\lVert u \rVert. \label{fnormdu}
\end{align} 
Using \cref{basisCond} and \cref{normE}, we can bound $\mE$ and $(z\mI -\tilde{\mL})|_{\Phi_n}$ as 
\begin{align*}
\lVert \mE \rVert \leq \frac{\beta^{\phi}_n}{\alpha^{\psi}_{n+b}} \gamma_{bn^{3/2}} \lVert A_{11} \rVert,~~~~ \lVert (z\mI -\tilde{\mL})|_{\Phi_n}\rVert \geq \frac{\alpha^{\phi}_n}{\beta_{n+b}^{\psi}}\lVert A_{11} \rVert.
\end{align*}
These two inequalities gives 
\begin{align}
\lVert \mE \rVert \leq \kappa^{\phi}_n\kappa^{\psi}_{n+b}\gamma_{bn^{3/2}}\lVert (z\mI -\tilde{\mL})|_{\Phi_n}\rVert. \label{normmE}
\end{align}
Since $\lVert \mR(z) \Delta \mL \rVert < 1/2$, we have
\begin{align*}
\lVert (z\mI -\tilde{\mL})|_{\Phi_n}\rVert \leq \lVert (z\mI -\mL)|_{\Phi_n}\rVert \lVert (\mI + \mR(z)\Delta\mL) \rVert \leq \frac{3}{2} \lVert (z\mI -\mL)|_{\Phi_n}\rVert.
\end{align*}
This enables us to replace $\lVert (z\mI -\tilde{\mL})|_{\Phi_n}\rVert$ in \cref{normmE} by $\lVert (z\mI -\mL)|_{\Phi_n}\rVert$ to have
\begin{align}
\lVert \mE \rVert \leq \kappa^{\phi}_n\kappa^{\psi}_{n+b}\gamma_{bn^{3/2}}\lVert (z\mI -\mL)|_{\Phi_n}\rVert,\label{fnormE}
\end{align}
where $3/2$ is absorbed into $\gamma_{bn^{3/2}}$. Subtracting \cref{opperteq} from \cref{resEqn} leads to
\begin{align}
(z\mI - \mL)(\mI + \mR(z) \Delta \mL + \mR(z)\mE)(v - \hat{v}) = \Delta \mL v + \mE v + \Delta u. \label{opresid}
\end{align}
Since $\lVert \mR(z) \Delta \mL \rVert < 1/2$, and $\lVert \mR(z) \rVert\lVert \mE \rVert < 1/2$, it follows that
\begin{align*}
\lVert \mR(z) \Delta \mL + \mR(z)\mE \rVert < 1,
\end{align*}
which suggests the invertibility of $\mI + \mR(z) \Delta \mL + \mR(z)\mE$. Thus, we rewrite \cref{opresid} as
\begin{align*}
v - \hat{v} = (\mI + \mR(z) \Delta \mL + \mR(z)\mE)^{-1}\mR(z)(\Delta \mL v + \mE v + \Delta u),
\end{align*}
which leads to
\begin{align*}
\lVert v - \hat{v} \rVert \leq \frac{\lVert \mR(z) \rVert (\lVert \mR(z) \Delta \mL \rVert \lVert u \rVert + \lVert \mR(z) \rVert \lVert \mE \rVert \lVert u \rVert + \lVert \Delta u \rVert) }{(1 - \lVert \mR(z) \Delta \mL \rVert -\lVert \mR(z) \rVert\lVert \mE \rVert)}, 
\end{align*}
where we have used $\lVert v \rVert \leq \lVert \mR(z) \rVert\lVert u \rVert$. The last inequality, along with the assumption $\lVert \mR(z) \Delta \mL \rVert < \lVert \mR(z) \rVert\gamma_1$, \cref{fnormdu}, and \cref{fnormE}, gives \cref{errorv}.
\end{proof}
Since $\xi$ is large, the bound \cref{errorv} suggests that the error in $\hat{v}$ is sizable due to the ill-conditioning of \cref{resEqn}.


To streamline the subsequent discussion, we denote by
\begin{align}
(z^*\mI - \tilde{\mL}^*) \tilde{w} = \hv \label{adjEqnt}
\end{align}
the approximate operator equation obtained by substituting the variable coefficients or kernel function in \cref{adjEqn} with their finite normalized Legendre series expansions, and by
\begin{align}
B \tilde{\mfw} = \hat{\mfv} \label{adjsystem}
\end{align}
the resulting infinite-dimensional system. In other words, \cref{adjEqnt,adjsystem} are the adjoint counterparts of \cref{resEqnt,system}, respectively. Again, $\tilde{\mL}^* \neq (\tilde{\mL})^*$. Additionally, let $\hbw$ and $\hat{w}$ be the computed solution to \cref{adjsystem} and the corresponding polynomial, respectively.

The next corollary follows from \cref{thm:v-vt} and the analogous bound for $\hat{w}$ by taking into account the error in $\hat{v}$ when it serves as the right-hand side of \cref{adjEqnt}. We omit the proof.
\begin{corollary}\label{thm:w-wt}
Suppose that $\bar{C}_1$, $\bar{C}_2$, and $\bar{C}_3$ are adjoint the counterparts of $C_1$, $C_2$, and $C_3$ in \cref{thm:v-vt}. The forward error in $\hat{w}$
\begin{align}
\lVert w - \hat{w} \rVert \leq  ((\bar{C}_1 + \bar{C}_3C_1)\xi^{3/2} + (\bar{C}_2 + \bar{C}_3C_2)\xi) \lVert u \rVert. \label{errw}
\end{align}
\end{corollary}

The key message conveyed by \cref{errw} can be expressed as 
\begin{align}
\lVert w - \hat{w} \rVert \le \zeta(u)\xi^{3/2}\epsilon_m\lVert u \rVert,\label{errorws}
\end{align}
where $\zeta(u)$ incorporates the non-essential factors. Here, we make the dependence on $u$ explicit for the lemma that follows.

Now we are in a position to give an error analysis for \cref{alg:operator}. There are three sources of error in \cref{alg:operator,alg:lanczos} that contribute to the error in the computed resolvent norm:
\begin{itemize}[leftmargin=*, itemsep=0pt]
\item The error in $\hat{w}$ incurred in solving \cref{Eqn}, quantified by \cref{errorws}.
\item The use of a finite-dimensional $H_k$, i.e., the stopping criterion in line 5 of \cref{alg:operator}.
\item The rounding errors occur and snowball elsewhere throughout the computation.
\end{itemize}

For convenience, we now denote by $e_S$, $e_T$, and $e_R$ these sources. Since the error in $\hat{w}$ is enormous, as suggested by \cref{errorws}, $e_S \gg e_R$. It is therefore safe to ignore $e_R$ as we do in the analysis below. Whether $e_T$ plays a role depends on when and how the Lanczos process is terminated. Thus, it is preferable to first examine how $e_S$, i.e., the large error in $\hat{w}$, is propagated throughout the Lanczos process before identifying an adaptive stopping criterion so that $e_T$ is comparable in magnitude with $e_S$. This allows for the earliest possible termination of the Lanczos process without compromising the accuracy of the computed resolvent norm.

In the remainder of this section, we concentrate only the numerical error in the computed version of $\mU_k$ and regard other quantities and the related computations as exact. This implies that we should have $\mU_k$ adorned with hat. However, we drop the hat from now on to keep the notations less busy. The lemma below can be deemed as an extension of Paige's theorem \cite{pai1} to operators acting on complex-valued functions. We therefore omit the proof. The first bound on the residual error $\Delta \mU_k$ is derived from \cref{errorws}, and it assesses the extent to which \cref{lanczos} fails to hold exactly. The second bound is deduced from the first, measuring the extent to which $\mU_k$ loses orthogonality. 
\begin{lemma}\label{lem:PerAna}
If $\zeta^{(k)} = \max(\zeta(u_1), \zeta(u_2),\dots, \zeta(u_k))$, 
\begin{align*}
\mT(z)\mU_k = \mU_kH_k+\beta_{k+1}u_{k+1}e_k^T +\Delta\mU_k , \label{MatErr}
\end{align*}
where $\Delta\mU_k  = (\Delta u_1\vert \dots\vert \Delta u_k)$ with $\Vert \Delta u_{k} \Vert \le \zeta^{(k)}\xi^{3/2}\epsilon_m$. Let $R_k \in \mathbb{C}^{k \times k}$ be the strictly upper triangular matrix with the $(i,j)$th entry $\rho_{ij} = \langle u_i, u_j \rangle$. Then
\begin{align*}
H_kR_k-R_kH_k = \beta_{k+1}\mU_k^*u_{k+1}e_k^{T} + \Delta R_k  \label{LossOrth}
\end{align*}
with $\lVert \Delta R_k  \rVert _F < \sqrt{2}k\zeta^{(k)}\xi^{3/2}\epsilon_m$.
\end{lemma}

Now we are in a position to show how accurately the eigenvalues of $H_k$ approximate those of $\mT(z)$ with the bounds for $\lVert \Delta \mU_k \rVert$ and $\lVert \Delta \mR_k \rVert_F$.
\begin{theorem}\label{thm:main}
Let the eigendecomposition of $H_k$ be
\begin{align*}
H_kY^{(k)} = Y^{(k)}diag(\mu_j^{(k)}),
\end{align*}
for which we assume $\mu_1^{(k)} > \mu_2^{(k)} > \cdots > \mu_k^{(k)}$. Let $y_j^{(k)}$ and $y^{(k)}_{ij}$ be the $j$th column and the $(i,j)$th element of the orthonormal matrix $Y^{(k)}$ respectively. We denote by $(\mu_j^{(k)}, z_j^{(k)})$ the corresponding Ritz pairs, where $z_j^{(k)} = \mU_k y_j^{(k)}$. If $\mu_{j}^{(k)}$ are well separated, i.e.,
\begin{align*}
\min_{i \neq j}|\mu_{j}^{(k)}-\mu_{i}^{(k)}| \geq \sqrt{2}k^{5/2}\zeta^{(k)}\xi^{3/2}\epsilon_m,
\end{align*}
the distance between $\mu_{j}^{(k)}$ and the true eigenvalue of $\mT(z)$ is then bounded as
\begin{align}
dist(\mu_{j}^{(k)}, \sigma(\mT(z))) < \frac{5}{2}\left(\beta_{k+1}|y_{kj}^{(k)}|+\sqrt{k}\zeta^{(k)}\xi^{3/2}\epsilon_m\right)\label{ErrorBound}
\end{align}
\end{theorem}

\begin{proof}
First, we note that the Ritz pair $(\mu_j^{(k)}, z_j^{(k)})$ satisfies
\begin{align}
\mT(z)z_j^{(k)} = \mu_j^{(k)}z_j^{(k)} + \beta_{k+1}y_{kj}^{(k)}u_{k+1} +\Delta\mU_k y_j^{(k)}. \label{pertlanczos}
\end{align}
Let $\mE_{j}^{(k)}:\mH \rightarrow \mH$ be the operator that effects $\mE_{j}^{(k)}x = -r_j^{(k)}\langle z_j^{(k)}, x \rangle/\langle z_j^{(k)}, z_j^{(k)} \rangle$. It follows from \cref{pertlanczos} that
\begin{align}
\left(\mT(z)+\mE_{j}^{(k)}\right) z_j^{(k)} = \mu_j^{(k)} z_j^{(k)}    \label{PretOp}
\end{align}
for $r_j^{(k)} = \beta_{k+1}y_{kj}^{(k)}u_{k+1} +\Delta\mU_k y_j^{(k)}$ and $\lVert r_j^{(k)} \rVert \le \beta_{k+1}|y_{kj}^{(k)}| + \sqrt{k}\zeta^{(k)}\xi^{3/2}\epsilon_m$.
Equation \cref{PretOp} shows that $\mu_j^{(k)} \in \sigma_{\lVert \mE_{j}^{(k)}\rVert}\left(\mT(z)\right)$, and since $\mT(z)$ is normal
\begin{align*}
dist(\mu_{j}^{(k)}, \sigma(\mT(z))) < \lVert \mE_{j}^{(k)} \rVert = \frac{\lVert r_j^{(k)} \rVert}{\lVert z_j^{(k)} \rVert} \le \frac{\beta_{k+1}|y_{kj}^{(k)}| + \sqrt{k}\zeta^{(k)}\xi^{3/2}\epsilon_m}{\lVert z_j^{(k)} \rVert}.
\end{align*}
\cref{lem:PerAna} and the results given in \cite[\S 3]{pai2} lead to \cref{ErrorBound}.
\end{proof}

\subsection{Stopping criterion}\label{sec:exit}
Suppose that $\mu_1^{(k)}$ converges to $\lambda_1(\mT(z))$, as described by \cref{thm:Convergence}. \cref{thm:main} supplies a bound for the relative error in $\mu_1^{(k)}$ as an approximation to $\lambda_1(\mT(z))$
\begin{align}
\frac{|\mu_1^{(k)}-\lambda_1(\mT(z))|}{|\lambda_1(\mT(z))|} \le \frac{5}{2}\left(\frac{\beta_{k+1}|y_{kj}^{(k)}|}{\xi}+\sqrt{k}\zeta^{(k)}\xi^{1/2}\epsilon_{m}\right), \label{relerr}
\end{align}
where we have used $\xi = \lambda_1(\mT(z))$ for $\mT(z)$ is a positive self-adjoint operator. 

Whereas the second term in the parentheses on the right-hand side of \cref{relerr} is the combination of $e_F$ and $e_S$ divided by $\lambda_1(\mT(z))$, i.e., the relative error owing to $\Delta \mL$, $E$, and $\Delta \mfu$, the first term is determined by the stopping criterion of the Lanczos process---given a tolerance $\delta$, we exit the Lanczos iteration when $\beta_{k+1}|y_{kj}^{(k)}|/\xi < \delta$. If $\delta > \sqrt{k}\zeta^{(k)}\xi^{1/2}\epsilon_{m}$, \cref{relerr} suggests that the relative error in the returned $\mu_1^{(k)}$ is $\mO(\delta)$. Otherwise, stopping after $\beta_{k+1}|y_{kj}^{(k)}|/\xi < \delta$ would incur unnecessary iterations. In this case, we should stop once $\beta_{k+1}|y_{kj}^{(k)}| < \sqrt{k}\zeta^{(k)}\xi^{3/2}\epsilon_{m}$ and the the relative error in $\mu_1^{(k)}$ is $\mO(\sqrt{k}\zeta^{(k)}\xi^{1/2}\epsilon_{m})$. Thus, we end up with the stopping criterion
\begin{align*}
{\beta_{k+1}|y_{kj}^{(k)}|} < \max(\sqrt{k}\zeta^{(k)}\epsilon_{m}\xi^{3/2},\ \delta\xi).
\end{align*}
Of course, there is no way that we can know $\xi$. Thus, we replace $\xi$ by $\mu_1^{(k)}$ to have
\begin{align}
{\beta_{k+1}|y_{kj}^{(k)}|} < \max\left(C_L\epsilon_{m}\left(\mu_1^{(k)}\right)^{3/2},\ \delta\mu_1^{(k)}\right), \label{sc}
\end{align}
where $C_{L}$ is a preset parameter that bounds $\sqrt{k}\zeta^{(k)}$. We set $C_{L}=100$ for all the experiments in \cref{sec:experiments} and find this value works uniformly well. It is not difficult to realized that \cref{sc} is analogous to the standard stopping criterion used in \textsc{Arpack} \cite[\S 4.6]{leh} and \textsc{Matlab}'s \texttt{eigs} function. In fact, \cref{sc} is virtually optimal in the sense that the best possible accuracy can be achieved if we set $\delta$ a value smaller than $\sqrt{k}\zeta^{(k)}\xi^{1/2}\epsilon_{m}$. See \cref{sec:laser} for an example demonstrating this.

When we are bailed out of the Lanczos process, the smallest possible relative error in $\mu_1^{(k)}$ is $\mO(\xi^{1/2}\epsilon_{m})$ if $\sqrt{k}\zeta^{(k)}$ is of modest magnitude. Thus, we can always have an idea on roughly how many faithful digits we have at most in $\sqrt{\mu_1^{(k)}}$ by its magnitude when the resolvent norm $\xi^{1/2}$ is no greater than $\mO(1/\epsilon_{m})$, i.e., $\xi^{-1/2} \gtrsim \mO(\epsilon_{m})$. If $\xi^{-1/2}$ turns out to be about $10^{-10}$ we then have roughly $3$ or $4$ digits to trust. If $\xi^{-1/2} \approx 10^{-20}$ is what we need, double precision is then off the table, and we have to look at quadruple or higher precisions. In such a case, we should repeat the entire computation with the extended precision, including the approximations to the variable coefficients and the kernel functions by normalized Legendre series. Usually, the quantities $m^{a_{N-1}},\ldots, m^{a_0}$, $m^K$, $m^f$, $m^g$, etc. are expected to become larger, since the tail of the series deemed to be negligible emerges (much) later in the extended precision arithmetic.

\section{Numerical experiments}\label{sec:experiments}
To demonstrate the advantages of the proposed method, we test it on a wide range of problems, covering all the operators we discussed in \cref{sec:operator}. The ``exact'' resolvent norm that we compare against in the remainder of this section is obtained using octuple precision\footnote{The $256$-bit octuple precision is the default format of \textsc{Julia}'s \texttt{BigFloat} type of floating point number. \texttt{BigFloat}, based on the GNU MPFR library, is the arbitrary precision floating point number type in \textsc{Julia}.} in \textsc{Julia}.

\subsection{Lasers}\label{sec:laser}
\begin{figure}[t!]
\centering
\subfloat[$\varepsilon$-pseudospectra of $\mF^L$]{\includegraphics[width=0.33\linewidth]{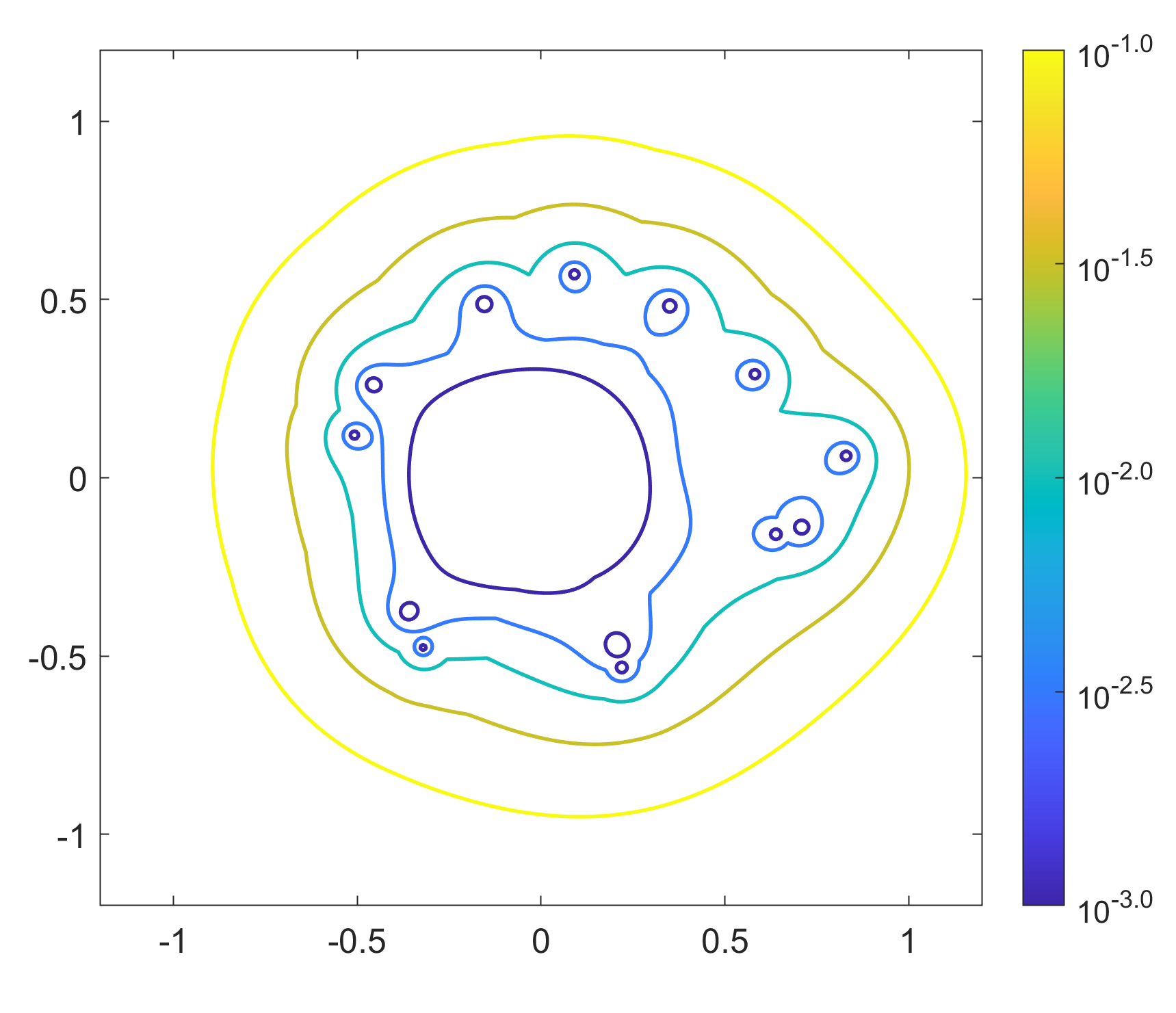}\label{fig:pseLaserF}}
\hfill 
\subfloat[$\varepsilon$-pseudospectra of $\mF_c^L$]{\includegraphics[width=0.33\linewidth]{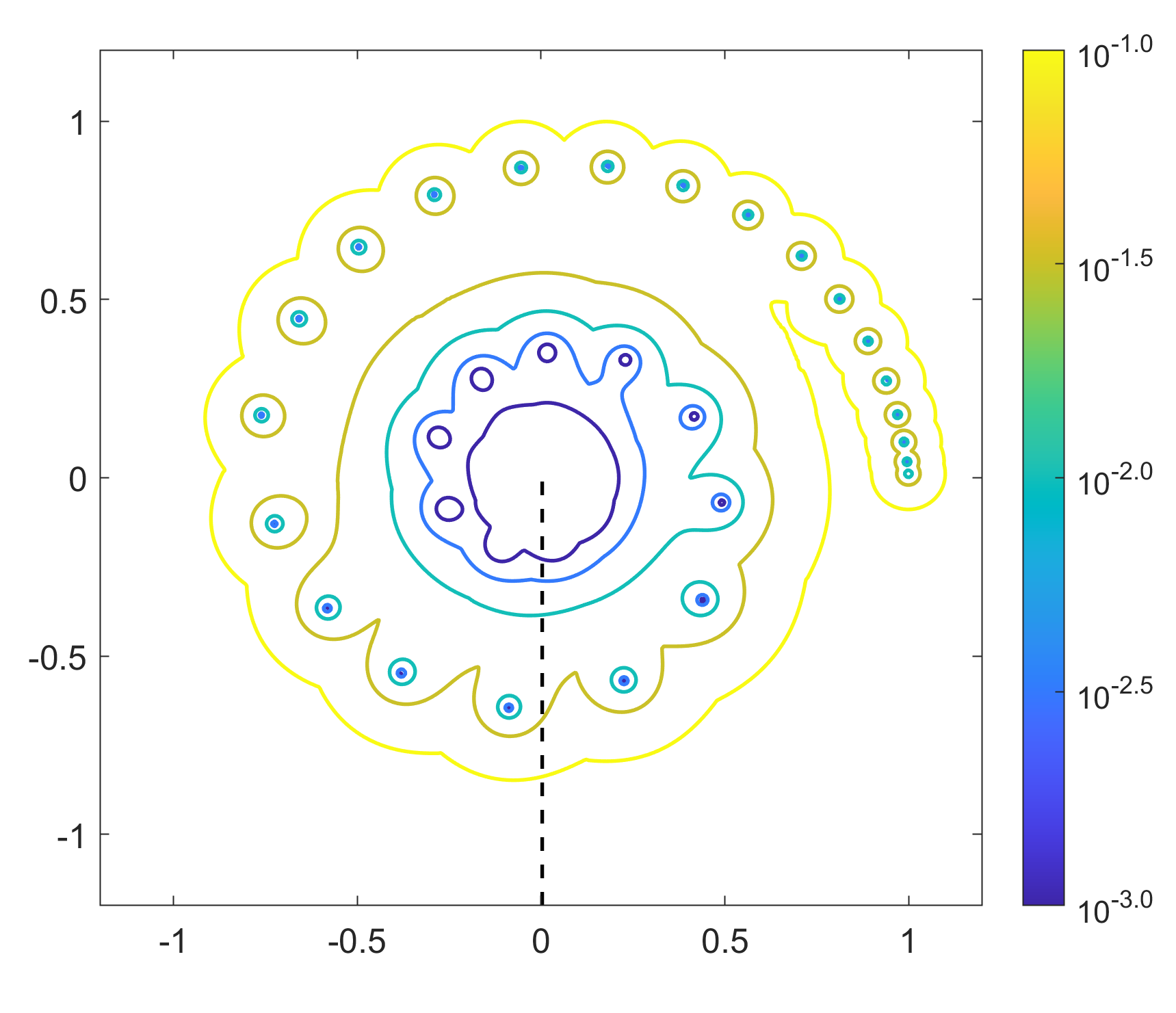}\label{fig:pseLaserC}}
\hfill 
\subfloat[relative error ]{\includegraphics[width=0.33\linewidth]{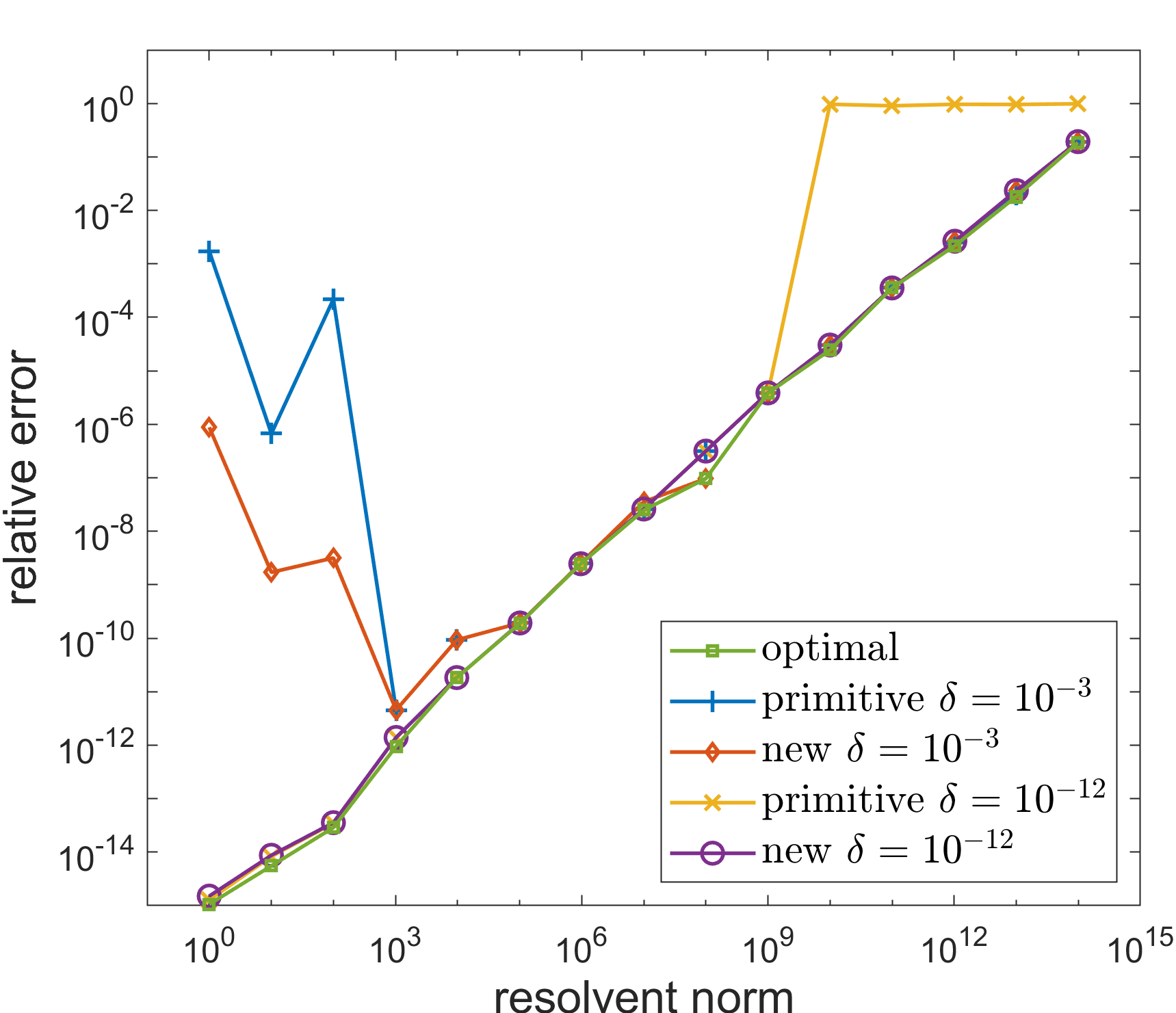}\label{fig:exitCondConv}}
\caption{Huygens--Fresnel operators}
\label{fig:pseLaser}
\end{figure}
Our first problems are the Huygens--Fresnel operators
\begin{align*}
(\mF^L u)(s) = \sqrt{\frac{iF}{\pi}}\int_{-1}^1  e^{-iFM(s/M-t)^2} u(t)\md t,~~~ (\mF_c^L u)(s) = \sqrt{\frac{iF}{\pi}}\int_{-1}^1  e^{-iF(s-t)^2} u(t)\md t
\end{align*}
from modeling the laser problem \cite[\S 60]{tre5}. These integral operators acting on $u \in L_2([-1,1])$ correspond to unstable and stable resonators respectively with the former being a general Fredholm operator and the latter a Fredholm operator of convolution type. In \cref{fig:pseLaserF,fig:pseLaserC}, we reproduce by \cref{alg:operator} the second panels of Figures 60.6 and 60.2 from \cite[\S 60]{tre5} for Fresnel number $F = 16\pi$ and magnification $M = 2$. The contours at $\varepsilon = 10^{-1}, 10^{-1.5}, \dots, 10^{-3}$ are a spot-on-match with the ones that Trefethen and Embree obtained using the traditional method.

\cref{fig:exitCondConv} shows the resolvent norm versus its relative error obtained with \cref{alg:operator} at the $15$ intersection points of the vertical dashed line $x = 0$ in \cref{fig:pseLaserC} and the boundaries of $\varepsilon$-pseudospectra for $\varepsilon = 10^{0}, 10^{-1}, 10^{-2}, \dots, 10^{-14}$ (not shown). The purple and red curves correspond to $\delta = 10^{-12}$ and $\delta = 10^{-3}$, respectively. We also search for the smallest relative error that can be achieved at each of these intersection points and have it plotted (green) for comparison. The green squares lie roughly on a straight line of slope $1$, confirming the analysis at the end of \cref{sec:error} that the best relative error of $\mu_1^{(k)}$ is $\mO(\xi^{1/2}\epsilon_m)$. It also can be seen that the deviation of the purple curve from the green one is tiny, suggesting that our adaptive stopping criterion is nearly optimal when $\delta$ is set to a sufficiently small value. The first $5$ diamonds on the red curve show early exits from the Lanczos process due to the low target accuracies, as expected. For comparison, we also include two curves obtained by replacing the adaptive stopping criterion \cref{sc} by the primitive one in line 6 of \cref{alg:core} for $\delta = 10^{-3}$ (blue) and $\delta = 10^{-12}$ (yellow) respectively. The blue curve behaves similarly like the red one---an early bail-out to save computation for a loose $\delta$. However, a $\delta$ that is excessively stringent like $10^{-12}$ results in non-convergence for $z$ with a large resolvent norm. This clearly shows that why the adaptive stopping criterion is indispensable for the success of \cref{alg:operator}.



\subsection{First derivative operator}\label{sec:first}
\begin{figure}[t!]
\centering
\subfloat[traditional method]{
\includegraphics[width=0.31\linewidth]{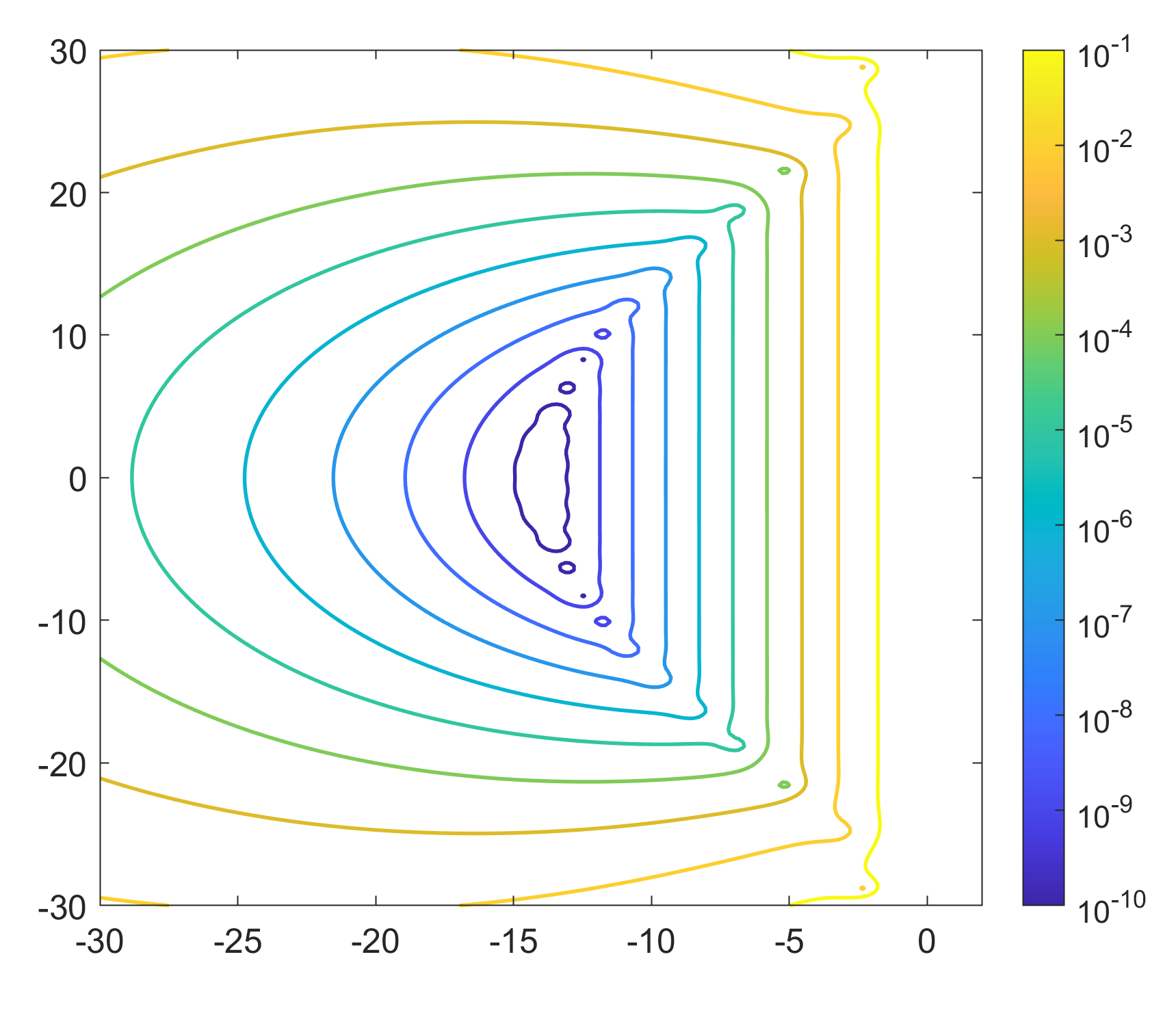}
\label{fig:pseFirstC}}
\hfill 
\subfloat[proposed method]{
\includegraphics[width=0.31\linewidth]{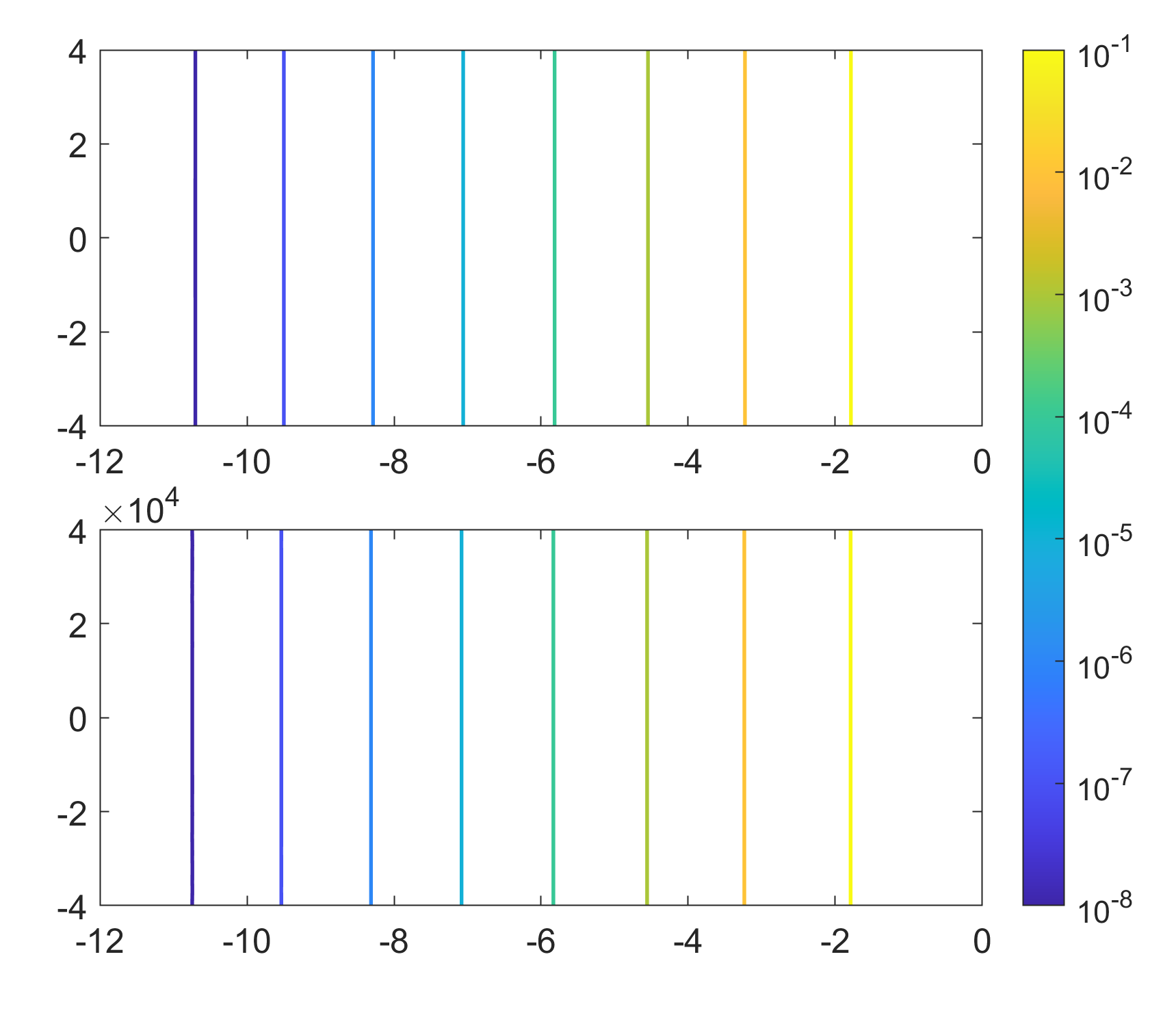}
\label{fig:pseFirstS}}
\hfill 
\subfloat[extending to the left]{
\includegraphics[width=0.31\linewidth]{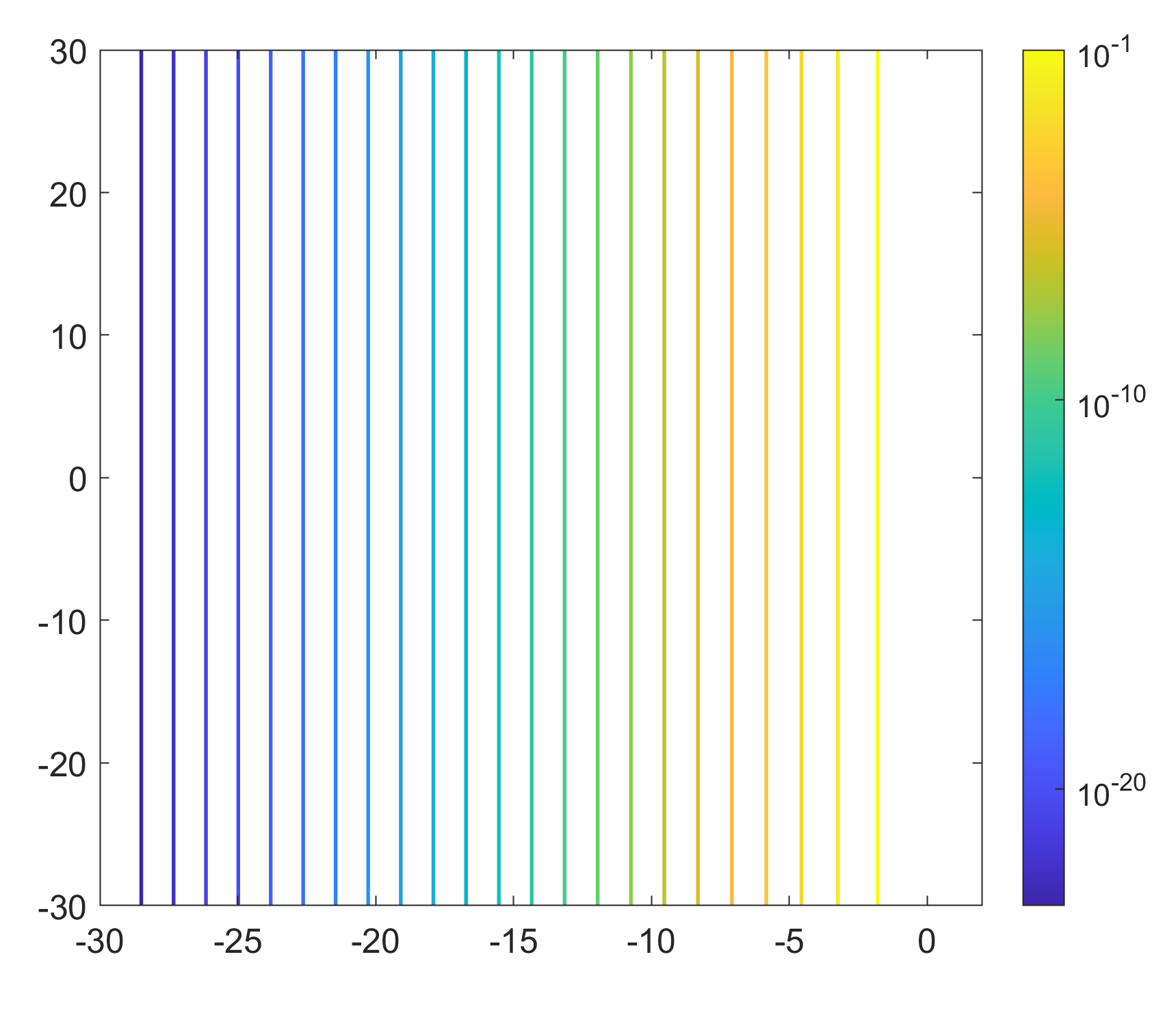}
\label{fig:pseFirstB}}
\caption{$\varepsilon$-pseudospectra of the first derivative operator.}
\label{fig:pseFirst}
\end{figure}
Our second example is the first derivative operator
\begin{align*}
\mQ^F u = \frac{\md u}{\md x},
\end{align*}
acting on $L^2([0,2])$ and subject to $u(2) = 0$ \cite{red2, tre2, tre5}. This simple yet interesting operator features nonexistence of eigenvalues and eigenmodes, and its resolvent norm depends only on $\Re(z)$, implying that the pseudospectra contours are straight vertical lines. The pseudospectra plot obtained by the traditional method is shown in \cref{fig:pseFirstC}. The discretize-then-solve paradigm is at least partially failed for this example, if not a fiasco---only the contours in a triangular region to the left of the origin resemble the true pseudospectra as vertical lines and the rest of the plot is an exhibition of spectral pollution. 

Whereas Figure 5.2 in \cite[\S 5]{tre5} and the upper panel of Figure 3 in \cite{tre2} are obtained from cropping a plot similar to \cref{fig:pseFirstC}, the upper panel of \cref{fig:pseFirstS} reproduces the same plot directly using the proposed method. In fact, we can go much beyond the $y$-scale in \cref{fig:pseFirstS}---the lower panel of \cref{fig:pseFirstS} also displays the same level sets, obtained by \cref{alg:operator}, but with the $y$-axis limit $10,000$ times greater. 

In the upper panel of \cref{fig:pseFirstS}, the left limit is set to $\Re(z) = -12$. We actually can push it a little further to $\Re(z)=-16.2$ (not shown) where the resolvent norm there is about $10^{12.5}$, the largest that we can compute with double precision arithmetic so that the two- or three-digit accuracy still allows the lines generated by the contour plotter to be straight to eyes. Since $\left\|(z\mI-\mQ^F)^{-1}\right\| = e^{2|\Re(z)|}/(2|\Re(z)|) + \mO \left(1/|\Re(z)|\right)$ for $\Re(z) < 0$, the resolvent norm grows virtually exponentially toward $-\infty$. As per the discussion of \cref{sec:exit}, if we intend to calculate the resolvent norm to the left of the vertical line corresponding to $\varepsilon = 10^{-12.5}$, we have to resort to higher precision arithmetic. \cref{fig:pseFirstB} shows the pseudospectra up to $x = -30$ and is meant to be compared with \cref{fig:pseFirstC}; the calculation is done in quadruple precision with $\epsilon_{m} \approx 1.93 \times 10^{-34}$.

\subsection{Advection-diffusion operator}\label{sec:advdiff}
\begin{figure}[t!]
\centering
\subfloat[$\varepsilon$-pseudospectra]{
\includegraphics[width=0.31\linewidth]{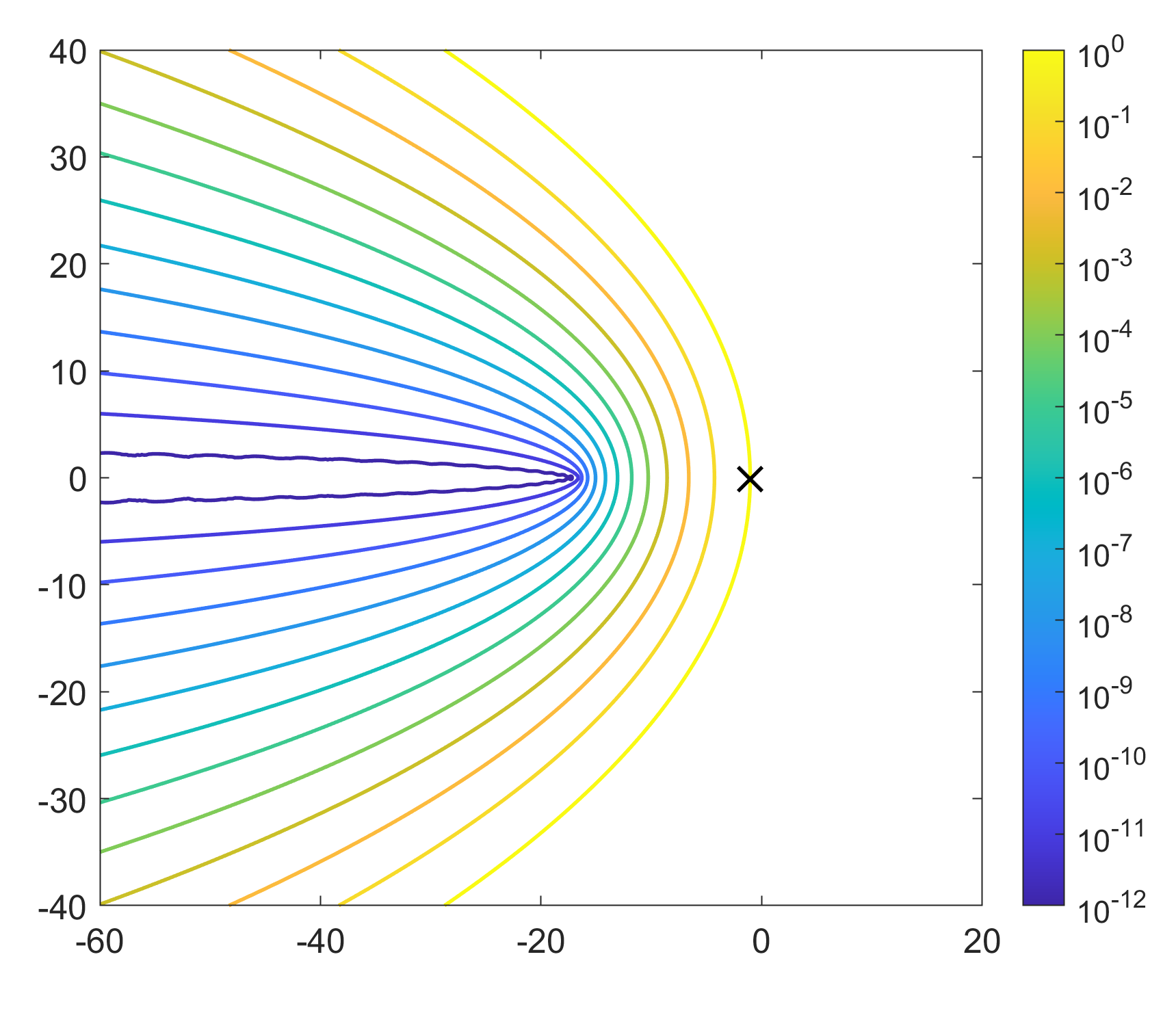}
\label{fig:pseConvL}}
\hfill 
\subfloat[spectral convergence]{
\includegraphics[width=0.31\linewidth]{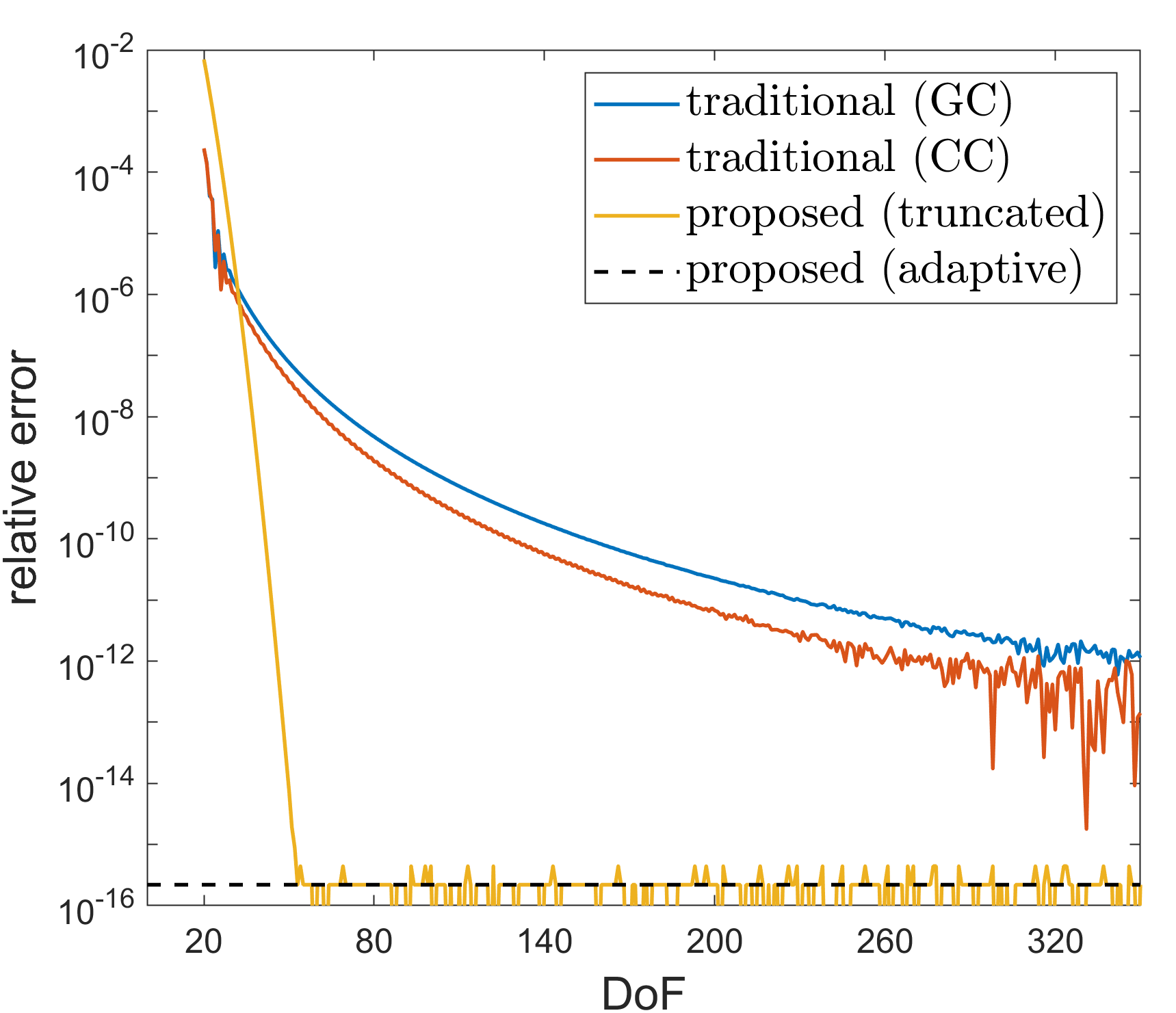}
\label{fig:relerrConv}}
\hfill 
\subfloat[close-up]{
\includegraphics[width=0.31\linewidth]{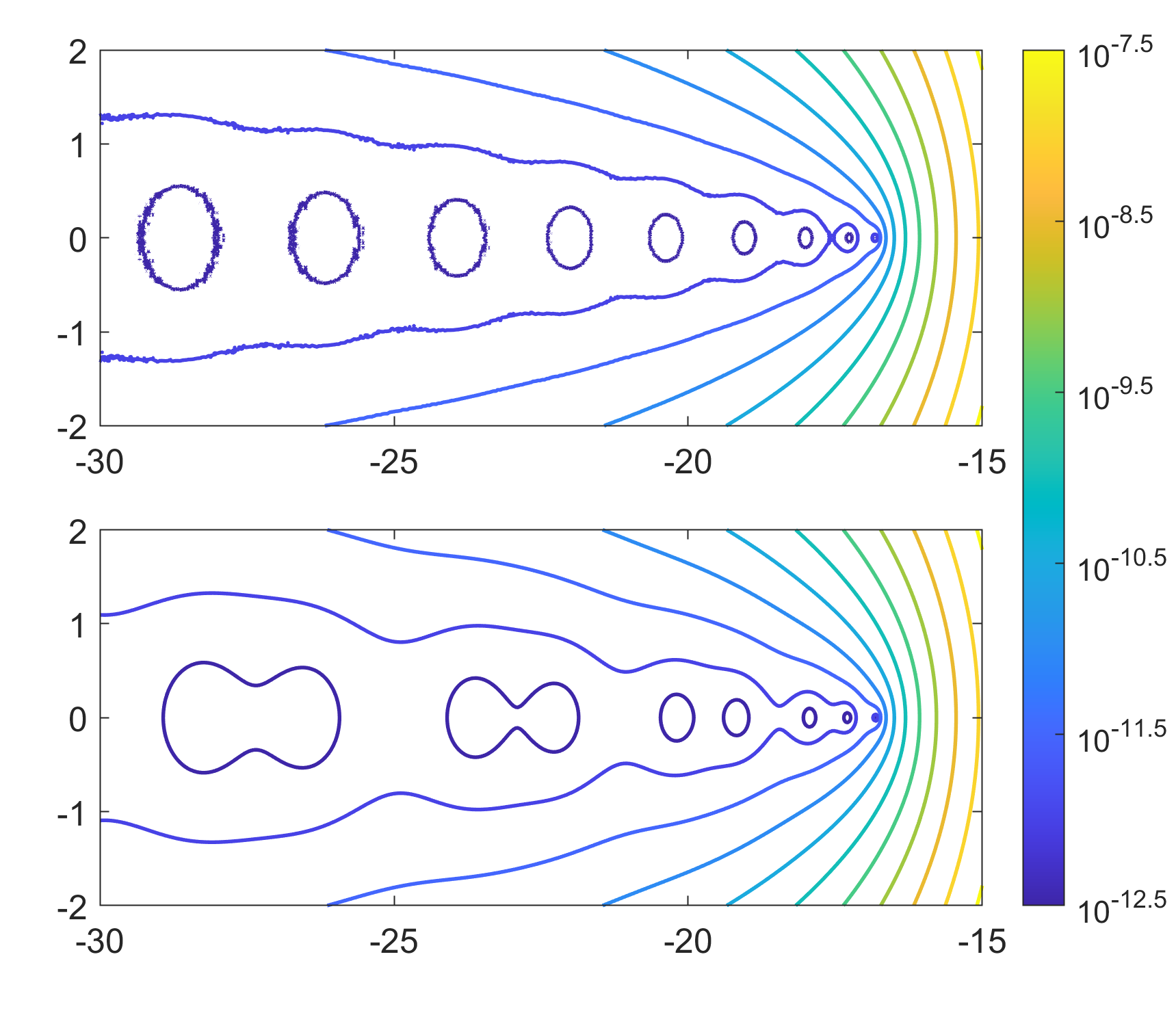}
\label{fig:pseConvM}}
\caption{Advection-diffusion operator}
\label{fig:pseConv}
\end{figure}

The advection-diffusion operator
\begin{align*}
\mQ^{AD} = \eta\frac{\md^2}{\md x^2} + \frac{\md}{\md x}
\end{align*}
acting on $L^2([0,1])$ and subject to the homogeneous Dirichlet boundary conditions serves in \cite{tre2} and \cite[\S 12]{tre5} as an example for transient effects caused by nonnormality. The pseudospectra of $\mQ^{AD}$ for $\eta = 0.015$ is shown as Figure 12.4 in \cite[\S 12]{tre5}, wheres \cref{fig:pseConvL} is a reproduction by \cref{alg:operator}. 

The advection-diffusion operator is also employed to demonstrate the quadrature issue related to the evaluation of the function norms \cite[\S 43]{tre5}. Trefethen and Embree suggest using a weight matrix corresponding to either Gauss--Chebyshev (GC) or Clenshaw--Curtis (CC) quadrature when Chebyshev pseudospectral method is used in a nonperiodic setting. However, they caution the loss of spectral convergence and a degeneration to algebraic convergence caused by the quadrature difficulty. This is confirmed by the blue and red curves in \cref{fig:relerrConv} which are the error-DoF plot of the resolvent norm at $-1.05-0.10i$ (cross in \cref{fig:pseConvL}). These curves are obtained using the traditional method with the GC and CC quadratures and very much resemble those given in Figure 43.5 of \cite[\S 43]{tre5} for the computed numerical abscissa. To the contrary, spectral convergence is gained by the proposed method using normalized Legendre polynomials. This can be seen from the yellow curve which is obtained by solving \cref{Eqn} with a fixed $n$ throughout the Lanczos iterations and plotting the corresponding error for various $n$---the error decays exponentially to machine precision at about $n=50$ then levels off, signaling adequate resolution. Of course, we do not have to choose $n$ manually in practice. The adaptivity in DoF (\cref{sec:solving}) helps us choose an optimal $n$ in each solution of \cref{Eqn} for complete resolutions. 
The resolvent norm at $-1.05-0.10i$ returned by our fully automated algorithm has an error no greater than $2.19\times 10^{-16}$, which is indicated by a dashed line in \cref{fig:relerrConv} for comparison.


\cref{fig:pseConvM} contains close-ups of \cref{fig:pseConvL} obtained using the proposed (upper) and the traditional (lower) methods respectively. Whereas the innermost contours in the upper panel barely hold their oval shape and start to wiggle, some of them coalesce in the lower panel which is totally erroneous. This is because the collocation-based matrix approximation to $\mQ^{AD}$ is ill-conditioned. It worsens the conditioning of the original problem and put the calculation at even higher risk. Thanks to the well-conditioned ultraspherical spectral methods, the precision is made the most of with only minimal loss of digits in solving \cref{Eqn}.

\subsection{GO and Wiener--Hopf operators}\label{sec:go}
\begin{figure}[t!]
\centering
\subfloat[\centering $\varepsilon$-pseudospectra of GO operator]{
\includegraphics[width=0.3\linewidth]{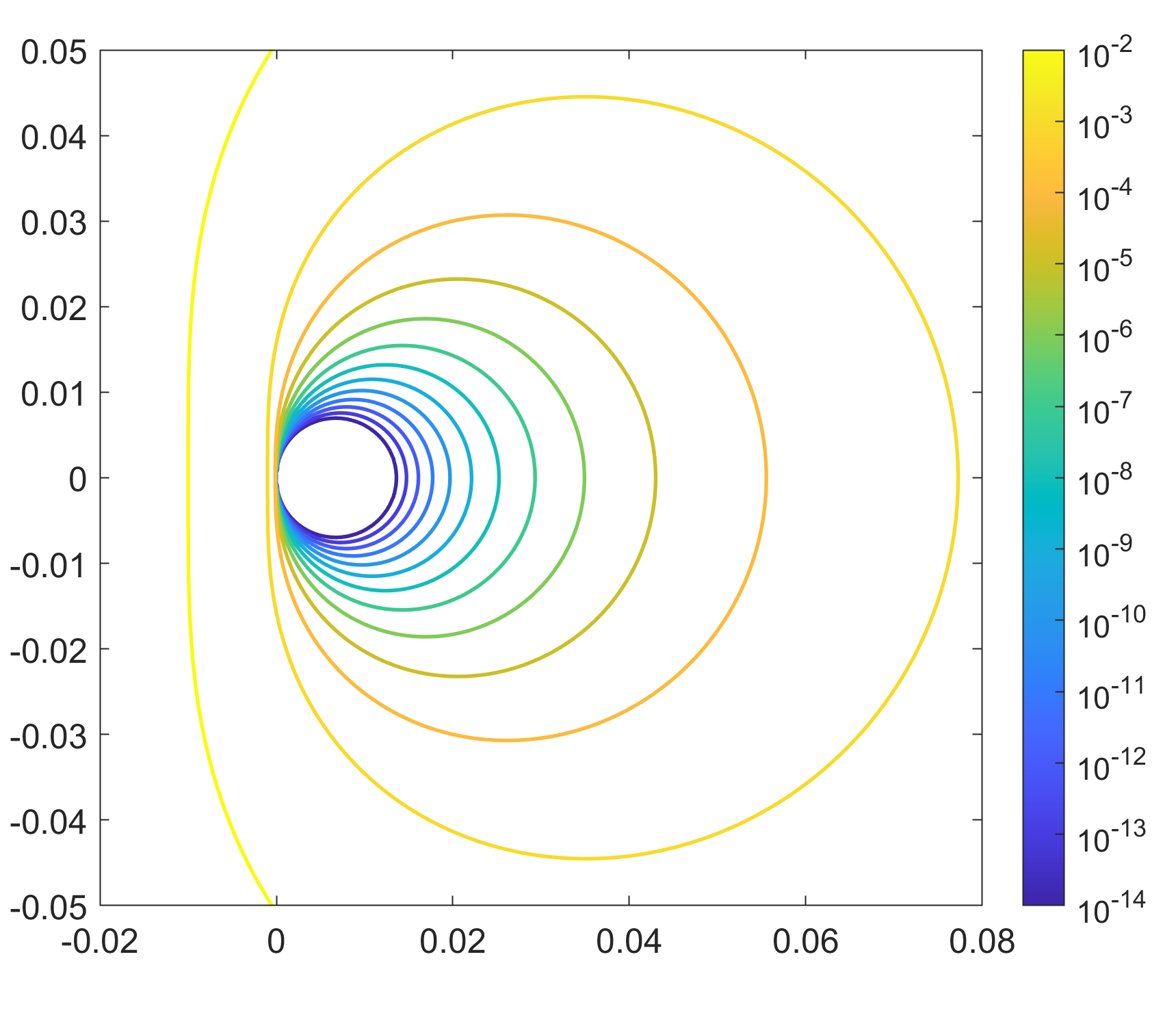}
\label{fig:pseOlver}}
\hfill 
\subfloat[\centering $\varepsilon$-pseudospectra of Wiener--Hopf operator]{
\includegraphics[width=0.3\linewidth]{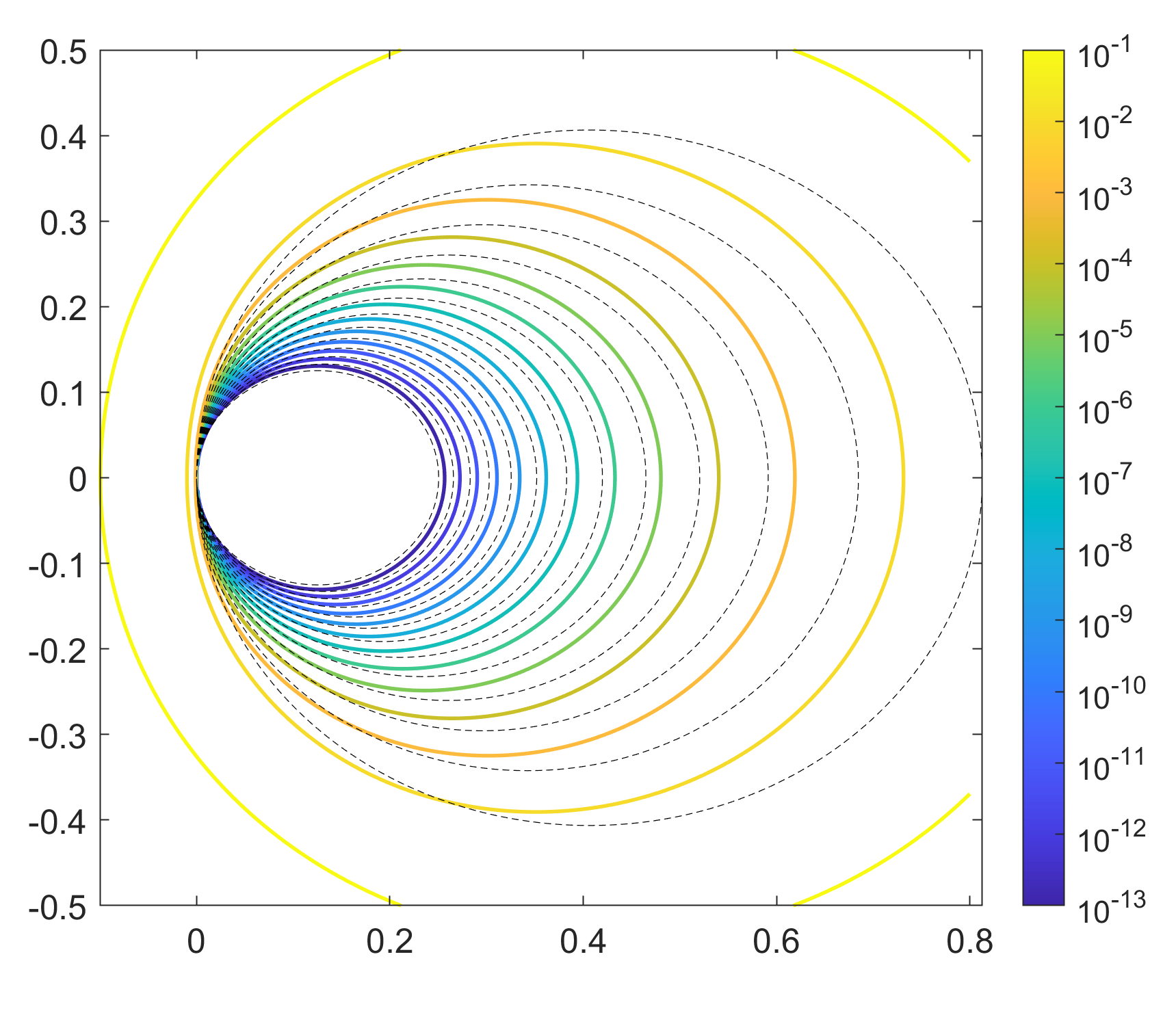}
\label{fig:pseWien}}
\hfill 
\subfloat[pointwise DoF]{
\includegraphics[width=0.3\linewidth]{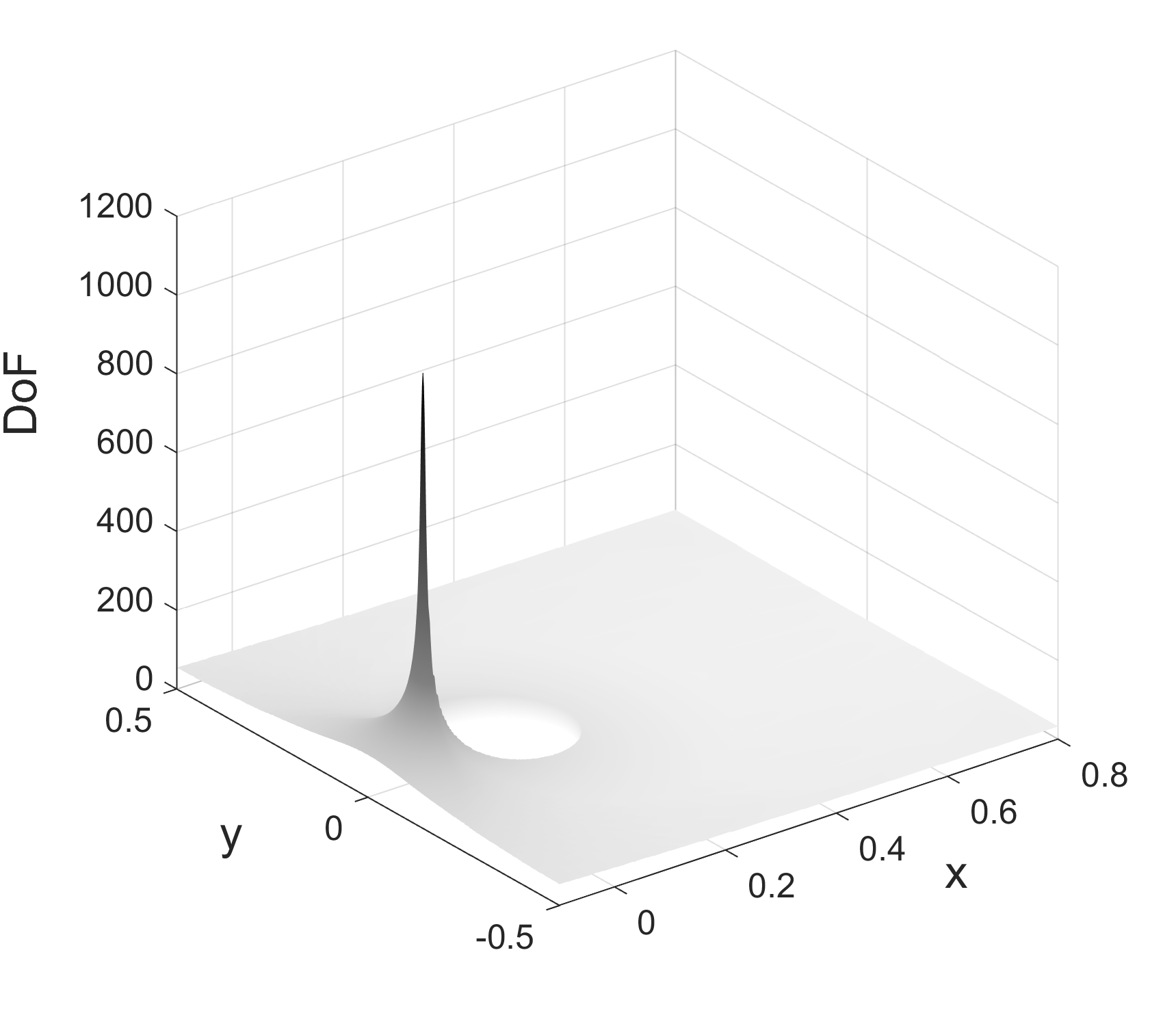}
\label{fig:dof}}\\
\subfloat[adaptive DoF]{
\includegraphics[width=0.3\linewidth]{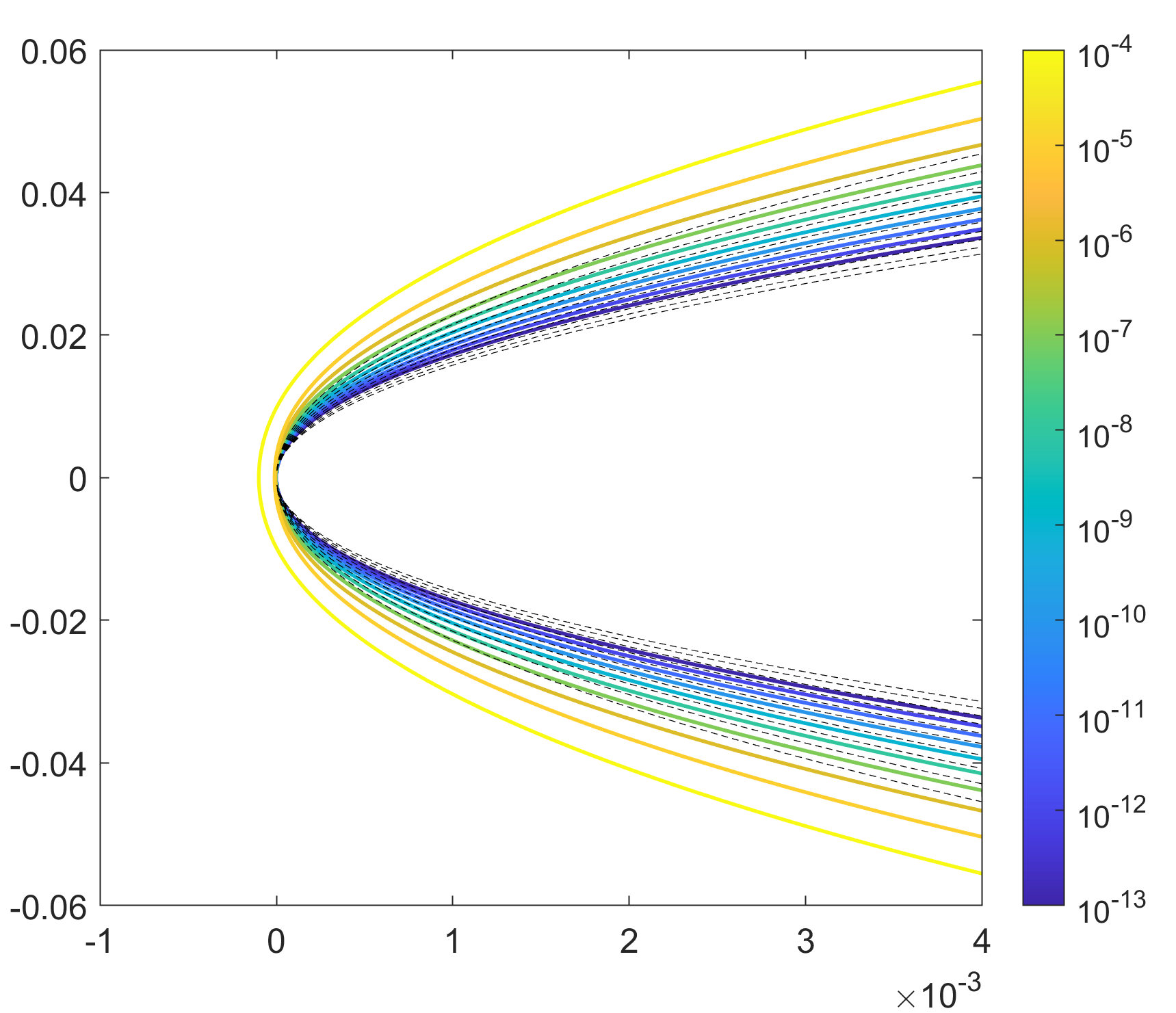}
\label{fig:pseWienT}}
\hfill 
\subfloat[pointwise DoF (close-up)]{
\includegraphics[width=0.3\linewidth]{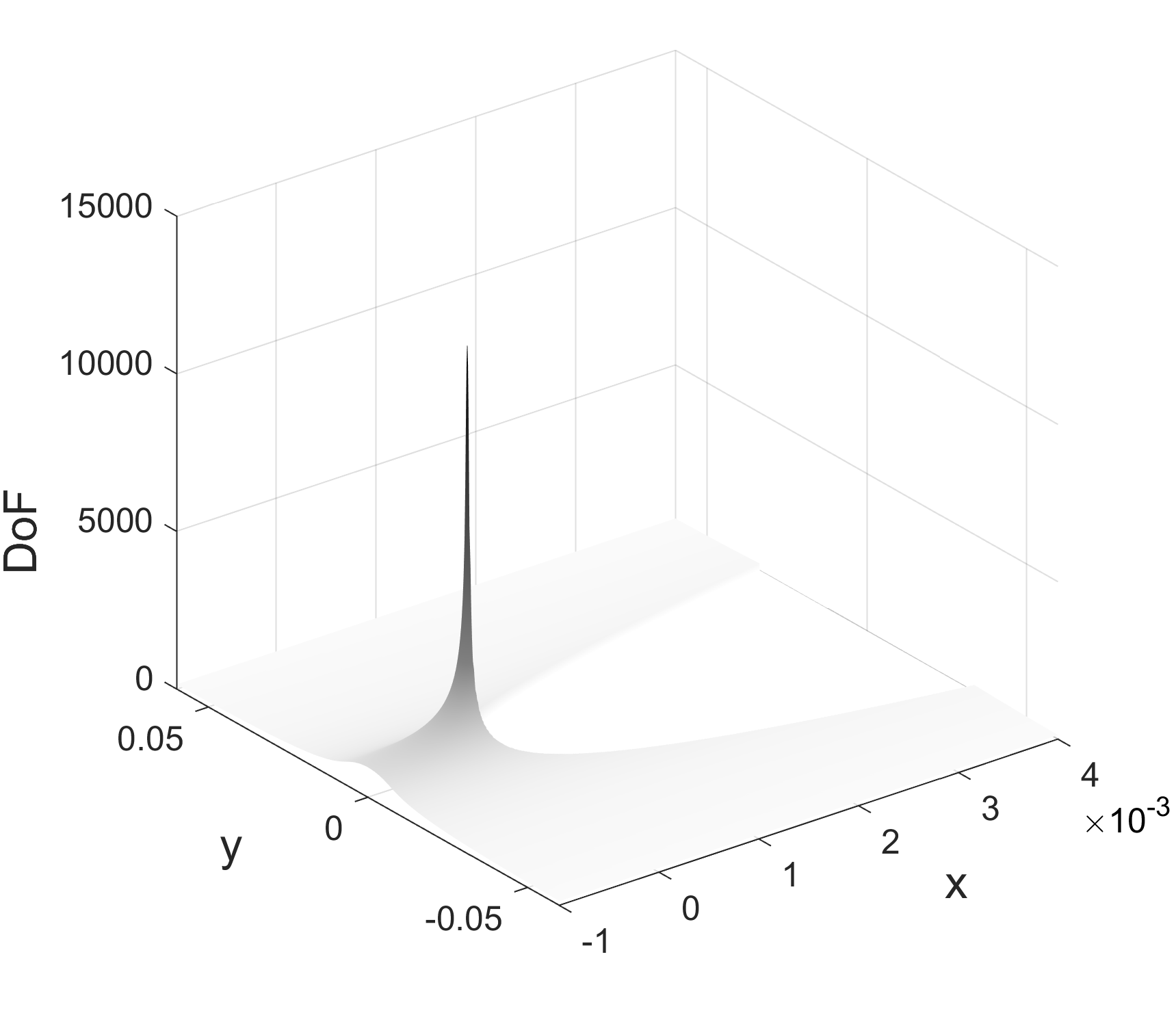}
\label{fig:dofS}}
\hfill 
\subfloat[$n = 200$]{
\includegraphics[width=0.3\linewidth]{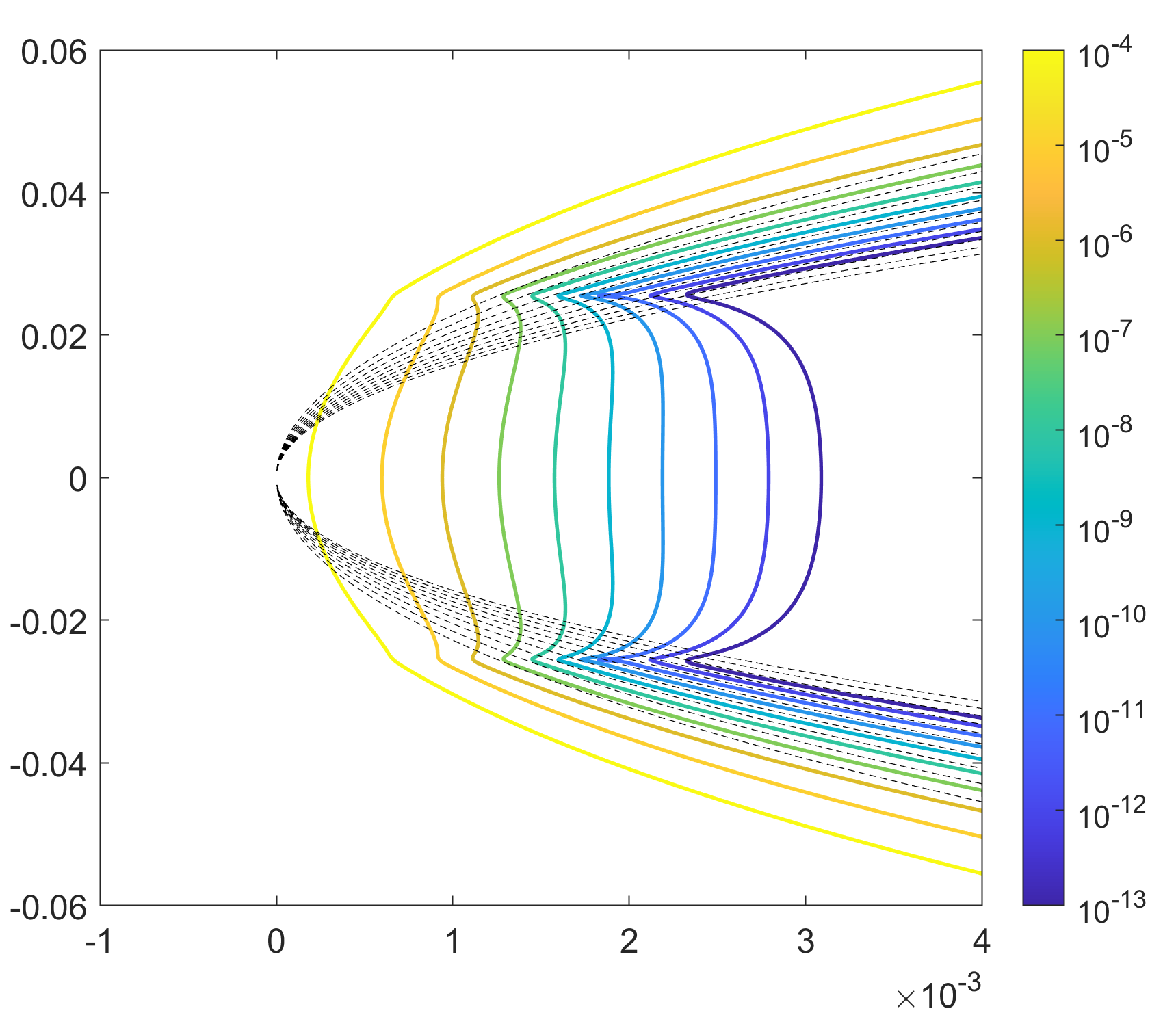}
\label{fig:pseWienF}}
\caption{GO and Wiener--Hopf operators}\label{fig:pseVolt}
\end{figure}

\cref{fig:pseOlver} shows the $\varepsilon$-pseudospectra of the GO operator \cite{gut}
\begin{align*}
(\mV^{GO} u)(s) = \int_{0}^{s}e^{-10(s-\frac{1}{3})^2-10(t-\frac{1}{3})^2}u(t)\md t, ~~~~ s \in [0, 1],
\end{align*}
which is computed via the approach elaborated in \cref{sec:voltimp}. Interestingly, we have not found any $\varepsilon$-pseudospectra plot of a general Volterra operator in the literature. This $\varepsilon$-pseudospectra plot resembles very much that of the Wiener--Hopf operator which is well-studied. The Wiener--Hopf operator
\begin{align*}
(\mV^{WH}_c u)(s) = \int_{s}^{d}e^{s-t}u(t)\md t, ~~~~ s \in [0, d]
\end{align*}
is a Volterra operator of convolution type. Previous study on its nonnormality and pseudospectra includes, e.g., \cite{red2, tre2}.

We show in \cref{fig:pseWien} the pseudospectra of the Wiener--Hopf operator $\mV_c$ for $d = 10$ by the boundaries of $\varepsilon$-pseudospectra for $\varepsilon = 10^{-1}, 10^{-2}, \dots, 10^{-13}$. These boundaries are in good agreement with their respective lower bounds (dashed curves) \cite{red2}. Specifically, the resolvent norm is calculated on a grid of $600 \times 600$ points before the results are sent to the contour plotter. For each grid point the largest $n$ used throughout the Lanczos process is shown in \cref{fig:dof}. It can be seen that $n$ grows as $z$ gets closer to the origin. At $-8.35 \times 10^{-4}\pm 8.35\times 10^{-4} i$, the closest grid points to the origin, the computed solution to \cref{Eqn} is expressed by a normalized Legendre series of degree $1031$, the highest among all the grid points. 

\cref{fig:pseWienT} is a close-up near the origin with the resolvent norm computed on a $800\times 800$ grid. As shown in \cref{fig:dofS}, $n = 13280$ for a complete resolution at $1.25 \times 10^{-6} \pm 7.51 \times 10^{-5}i$, the grid points nearest to the origin. Thanks to the adaptivity, the determination of $n$ is fully automated and sufficient resolution is guaranteed. In principle, the proposed method can compute the resolvent norm at a point arbitrarily close to the spectra, provided that sufficiently powerful hardware is available.

For comparison, \cref{fig:pseWienF} shows the boundaries of $\varepsilon$-pseudospectra obtained with $n=200$ throughout the entire grid. This ``finite section'' type implementation produces contours that exhibit conspicuous violation to the lower bounds for all the $\varepsilon$ considered, highlighting the importance of the adaptivity in DoF $n$.

\subsection{Orr--Sommerfeld}\label{sec:orr}
\begin{figure}[t!]
\centering
\subfloat[$L^2$ norm]{
\includegraphics[width=0.45\linewidth]{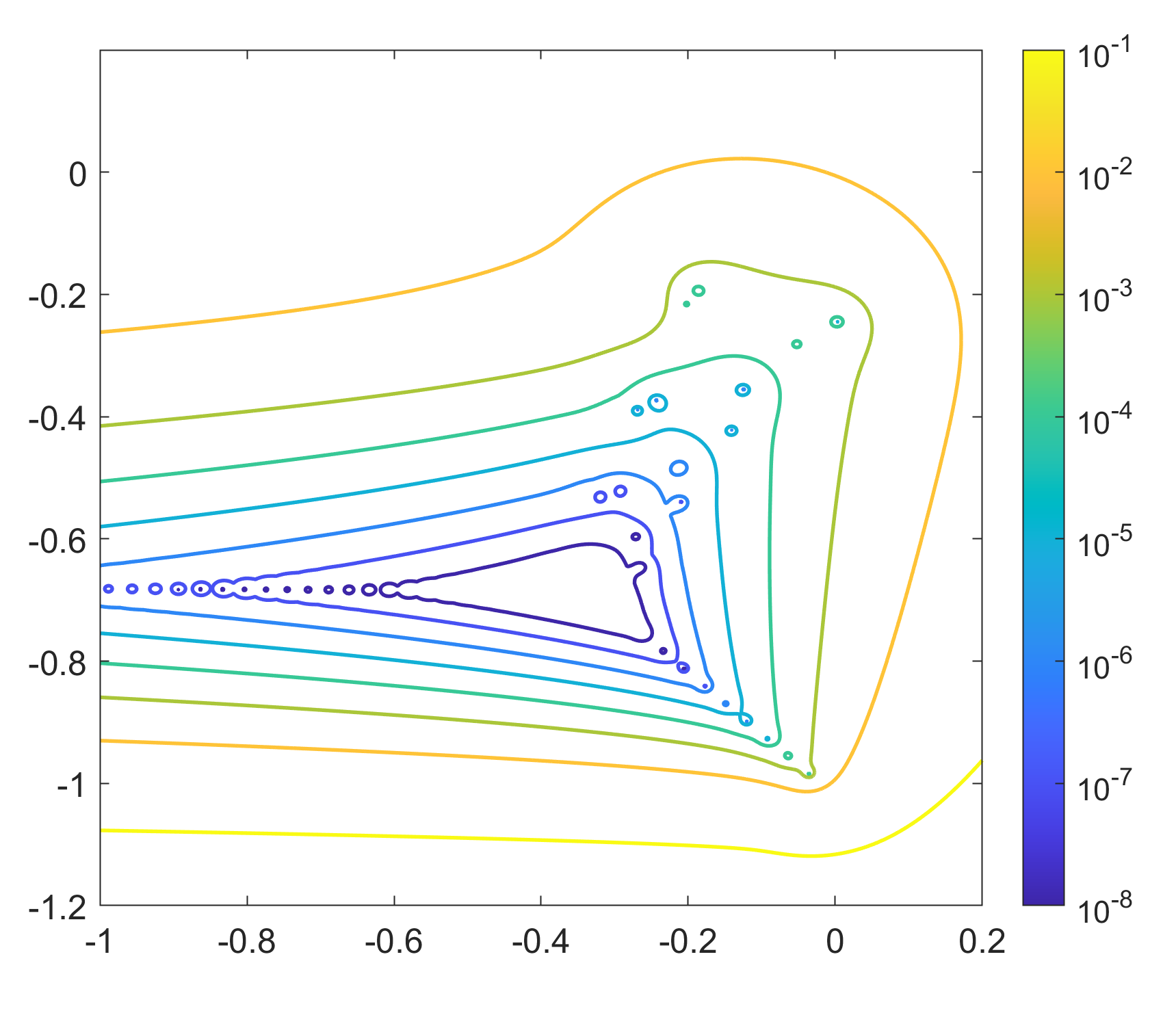}
\label{fig:pseOrrL2}}
\hfill 
\subfloat[energy norm]{
\includegraphics[width=0.45\linewidth]{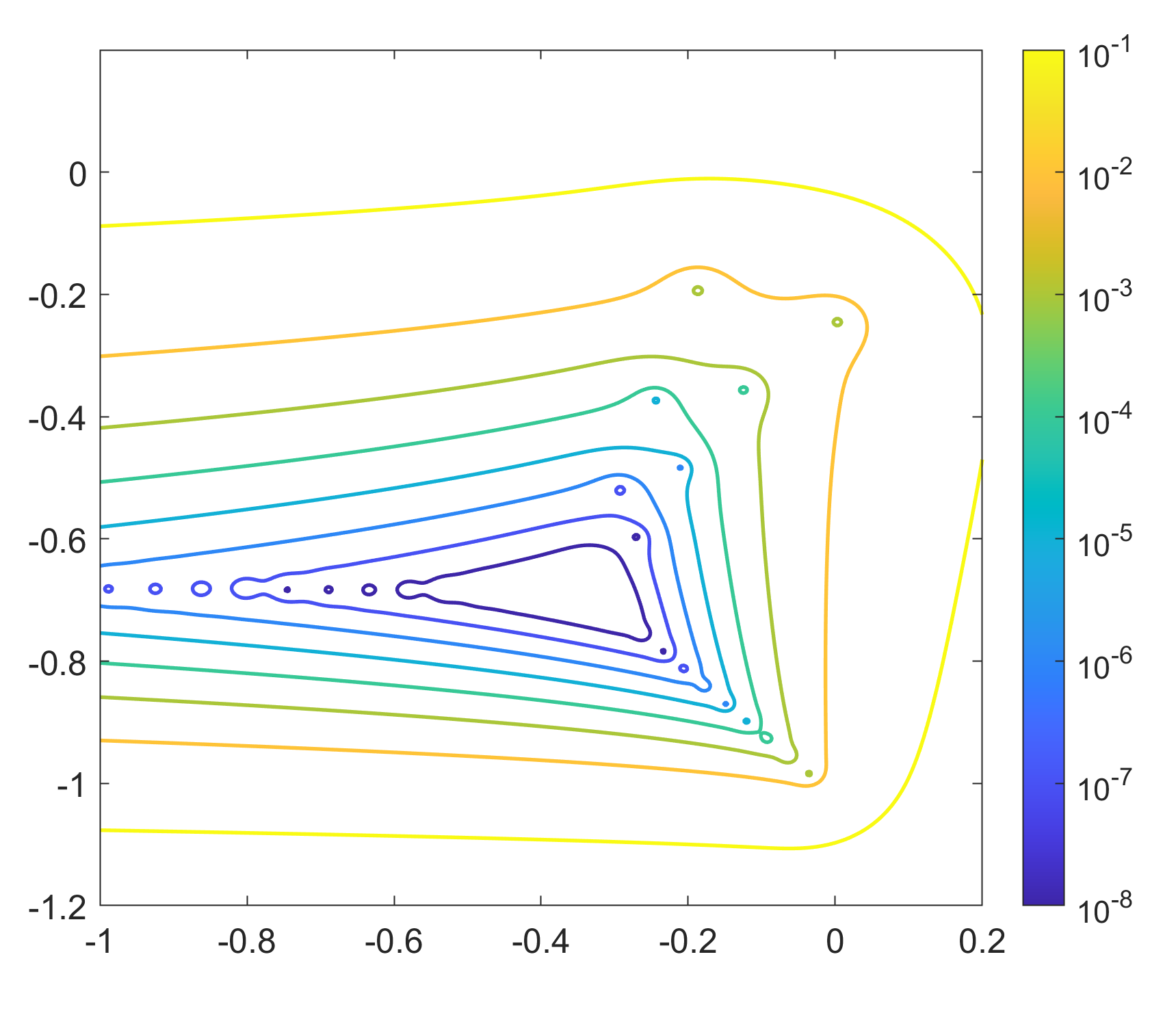}
\label{fig:pseOrrE}}
\caption{$\varepsilon$-pseudospectra of the Orr--Sommerfeld operator.}\label{fig:pseOrr}
\end{figure}
We close this section with the Orr--Sommerfeld operator for which pseudospectra plays a pivotal role in analyzing the temporal stability of fluid flows, as detailed in \cite{red3} and \cite[\S 22]{tre5}. In the standard form of \cref{gep}, the differential expressions of $\mA$ and $\mB$ read
\begin{align*}
\tau_{\mA} = -\frac1{R}\left(\frac{\md^2}{\md x^2}-\alpha^2\right)^2+i\alpha(1-x^2)\left(\frac{\md^2}{\md x^2}+\alpha^2\right)+2i\alpha, ~~~
\tau_{\mB} = -\frac{\md^2}{\md x^2}+\alpha^2
\end{align*}
with $\mA$ attached with $u(\pm 1)=0$ and $u'(\pm 1)=0$ and $\mB$ with $u(\pm 1)=0$. 

\cref{fig:pseOrrL2} displays the $2$-norm pseudospectra for $R = 10000$ and $\alpha = 1.02$ which is obtained following the details given in \cref{sec:gepimp}. If it is the pseudospectra in the energy norm \cite{red3} defined by
\begin{align*}
\langle \phi, \psi \rangle_E = \langle\mB\phi, \psi\rangle
\end{align*}
that we desire, we should note that $\mR^*(z)$ is different from \cref{R*} as \cref{R*} is obtained via the Euclidean inner product. Suppose that $\phi, \psi \in \mD(\mB)$. Reworking out as we do in \cref{inner} and noting $\mB = \mB^*$ for this particular example, we have 
\begin{align*}
\langle \mR(z)\phi, \psi \rangle_E = \langle \mB\phi, (z^*\mB-\mA^*)^{-1}\mB\psi \rangle = \langle \phi, (z^*\mB-\mA^*)^{-1}\mB\phi \rangle_E.
\end{align*}
Thus, $\mR^*(z) = (z^*\mB-\mA^*)^{-1}\mB$. Since $\mR(z)$ is not changed with the norm, \cref{R} still holds and we apply $\mR(z)$ by following \cref{Ru}. Then application of $\mR^*(z)$, i.e., $w = (z^*\mB-\mA^*)^{-1}\mB v$, is done by solving $(z^*\mB-\mA^*)w = \mB v$. \cref{fig:pseOrrE} is a reproduction of Figure 2 in \cite{red3}, showing the $\varepsilon$-pseudospectra in energy norm for $\varepsilon = 10^{-1}, 10^{-2}, \dots, 10^{-8}$.




\section{Extensions}\label{sec:extensions}
In this section, we discuss a few extensions of the proposed method, 
which, along with those dealt with in the previous sections, cover a great portion of the pseudospectra problems that one may come across in practice.

\begin{figure}[t!]
\centering
\subfloat[\centering 2D advection-diffusion operator]{
\includegraphics[width=0.32\linewidth]{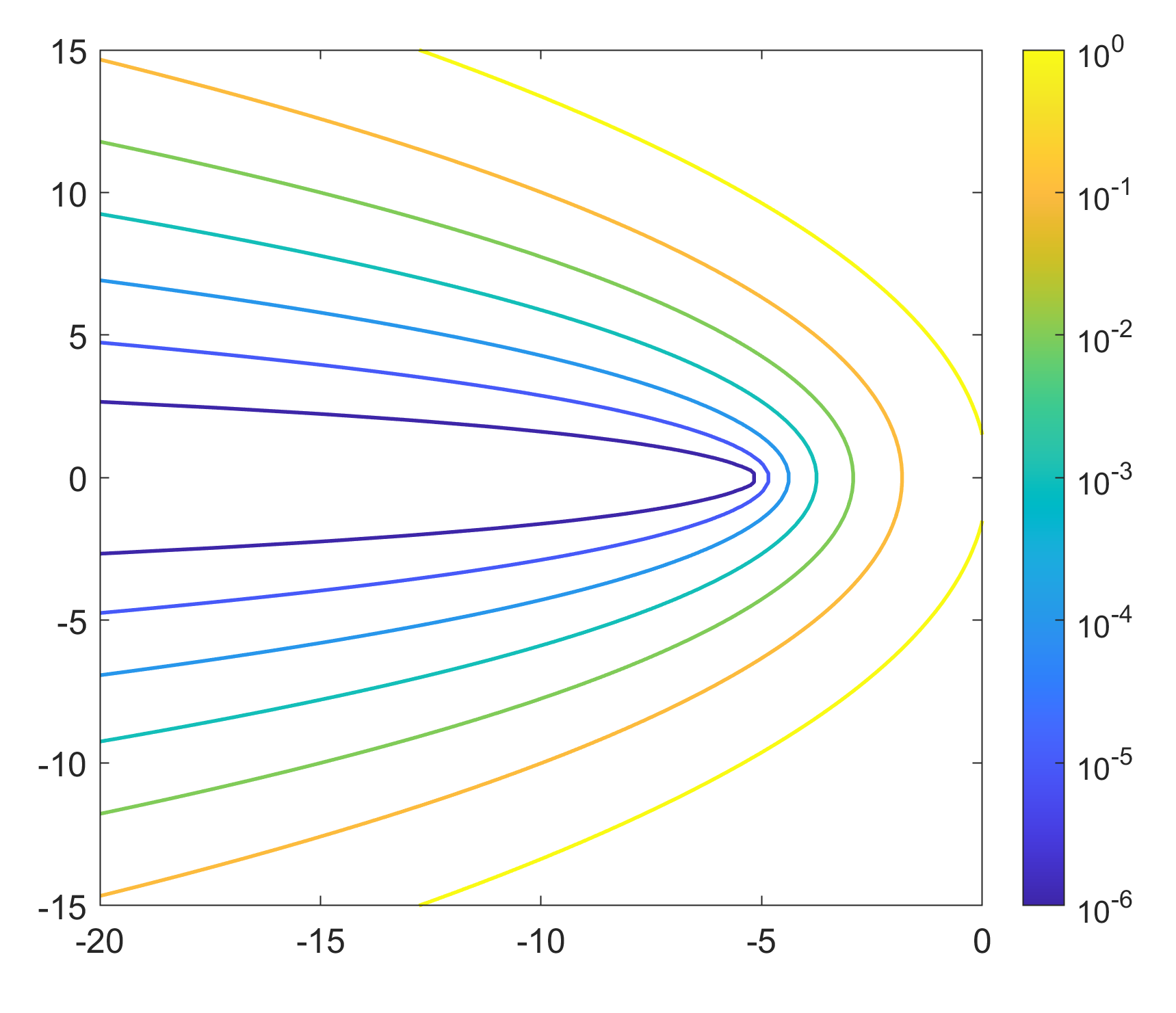}
\label{fig:pseConv2D}}
\hfill 
\subfloat[\centering block-structured wave operator]{
\includegraphics[width=0.32\linewidth]{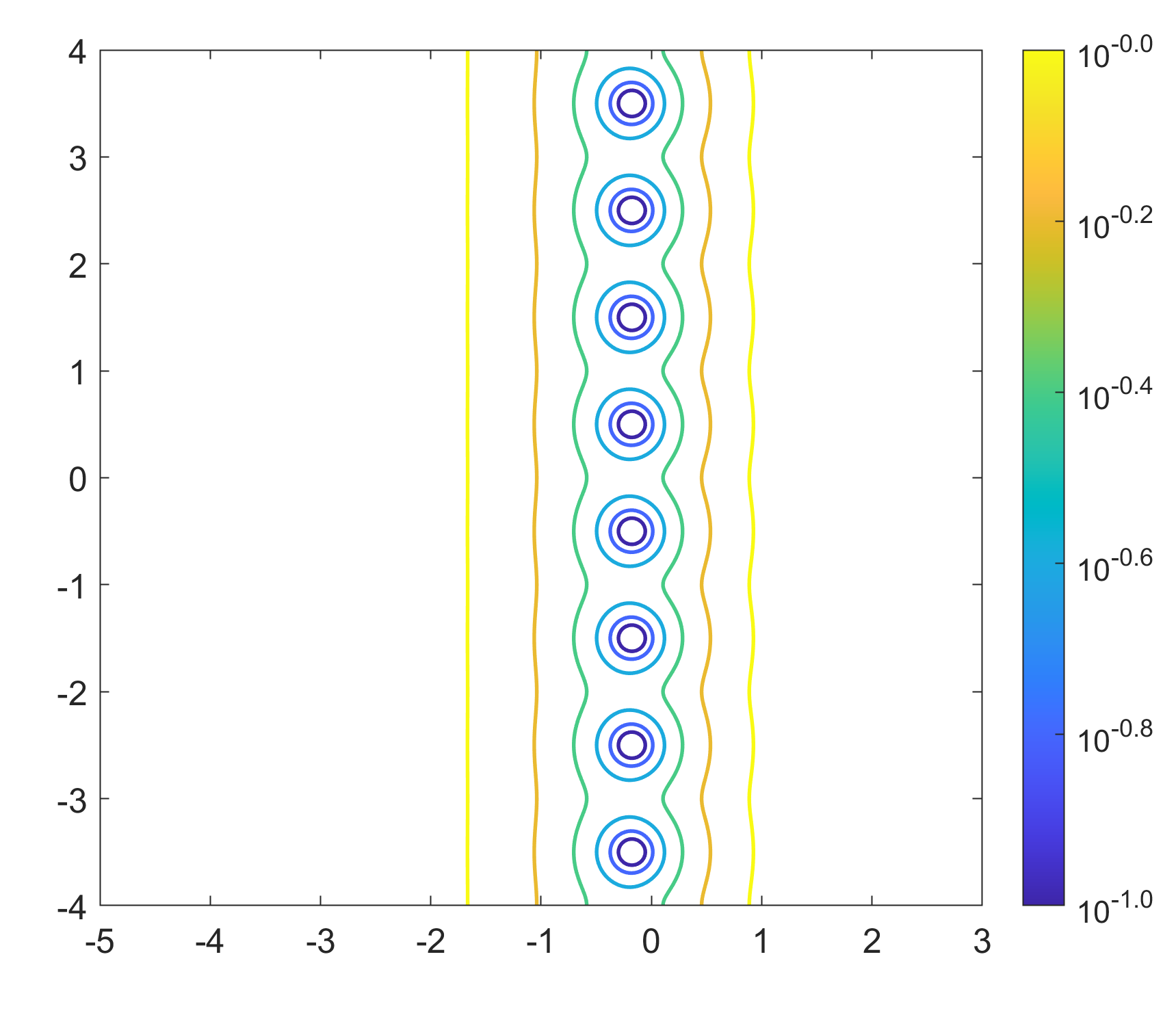}
\label{fig:pseWave}}
\hfill 
\subfloat[\centering Davies' operator in an unbounded domain]{\includegraphics[width=0.32\linewidth]{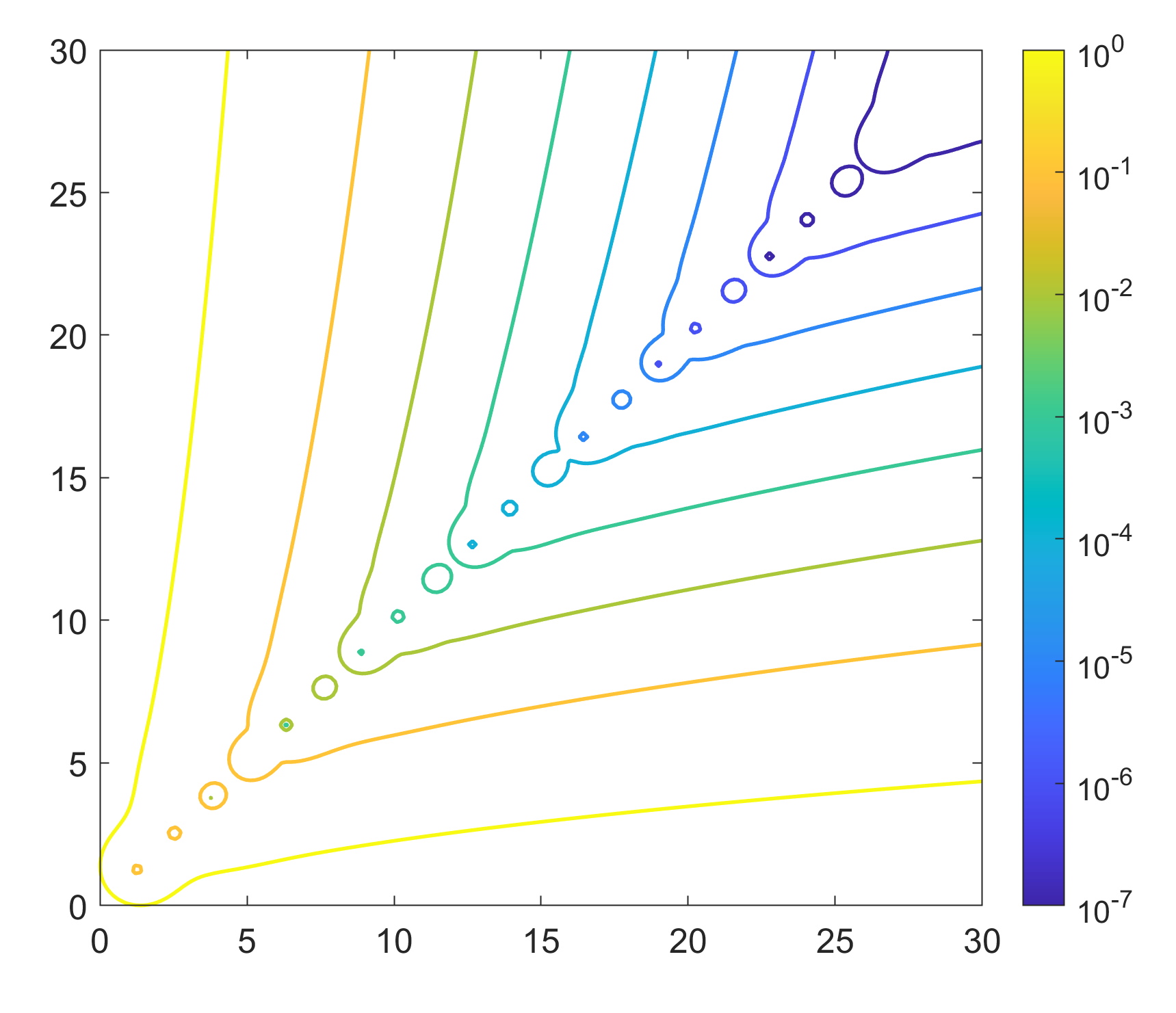}\label{fig:pseDavies2D}}
\caption{$\varepsilon$-pseudospectra of multivariate operator, block operator, and operator on unbounded domain.}
\end{figure}

\subsection{Multivariate operators}\label{sec:higher}
The existing strategies for handling an operator $\mL$ acting in multiple space dimensions include decoupling the pseudospectra problem into problems of lower dimensions \cite[\S 43]{tre5} and representing $\mL$ as an operator acting on $\ell^2(\mathbb{N})$ \cite{col3, col1}. On the other hand, the framework presented in this paper applies perfectly to an operator in multiple dimensions, provided that its resolvent is compact or compact-plus-scalar. We now take as an example the Dirichlet operator \cite[\S XIV.6]{goh2}, noting that most of the elliptic operators fall into this category. 


Given a bounded open set $\Omega$ in $\mathbb{R}^d$ and a partial differential expression
\begin{align*}
\tau_e = \sum_{i,j=1}^d a_{ij}(\mathbf{x}) \frac{\partial}{\partial x_i}\frac{\partial}{\partial x_j} + \sum_{i=1}^d b_i(\mathbf{x})\frac{\partial}{\partial x_i} + c(\mathbf{x})u,
\end{align*}
where each $a_{ij}$, $b_i$, and $c$ are in $C^{\infty}(\bar{\Omega})$ and $\bar{\Omega}$ is the closure of $\Omega$. The Dirichlet operator $\mP: L^2(\Omega) \rightarrow L^2(\Omega)$ is defined as $\mP u = \tau_e u$ with $\mD(\mP) = H_0^1(\Omega)\cap H^2(\Omega)$. If $\mP$ is uniformly elliptic, $\mR(z) = (z\mI-\mP)^{-1}$ is compact for $z \in \mathbb{C}\backslash \sigma(\mP)$ \cite[\S XIV.6]{goh2}.


In implementing \cref{alg:operator,alg:lanczos}, we solve \cref{Eqn} by the multivariate ultraspherical spectral method \cite{str, tow2}. The basis functions for such a problem are the tensor product of recombined normalized Legendre polynomials that satisfy boundary conditions in each spatial dimension. After reordering these tensor-product bases by the total degree, we obtain an infinite-dimensional matrix representation of $\mP$ with a growing bandwidth, which we solve by adaptive QR. \cref{fig:pseConv2D} shows the pseudospectra of the 2D advection-diffusion operator \cite{hem1}
\begin{align*}
\mP^{AD} u = \eta\Delta u+ \mathbf{V} \cdot \nabla u
\end{align*}
in $\Omega = [-1,1] \times [-1,1]$ with $\eta = 0.05$ and $\mathbf{V} = [0, -1]^{T}$.

\subsection{Block operators}\label{sec:OPsystem}
It is not uncommon to consider the pseudospectra of block operators that act on Cartesian product of function spaces. Such an operator may come from reducing a high-order problem to a system of low-order ones. For example, Trefethen and Driscoll in \cite{dri1,tre2} investigated the wave operator
\begin{align*}
\mQ^B = \begin{pmatrix}
 0 & \md/ \md x\\
 \md/ \md x & 0 
\end{pmatrix}
\end{align*}
with $\mD(\mQ^B) = \left\{(u, v)\in L^2([0, \pi]) \times L^2([0, \pi]) \mid v(0) = 0,\ u(\pi)+\eta v(\pi) = 0\right\}$ and $\eta = 1/2$. 
The compactness of $(z\mI-\mQ^B)^{-1}$ can be established following the proof in \cite[\S XI.1]{goh2}. \cref{fig:pseWave} displays the pseudospectra of $\mQ^B$ obtained using the proposed method. This replicates Figure 3 in \cite{dri1}, which was originally produced from a finite-dimensional discretization of $\mQ^B$ using the finite difference method. In solving \cref{Eqn} with the basis-recombined ultraspherical spectral method, each $\md/ \md x$ produces an infinite-dimensional banded matrix, thus representing $\mQ^B$ by a $2 \times 2$ block matrix of infinite dimensions with the $(1,2)$ and $(2,1)$ blocks banded. Applying a simple reordering gives an infinite-dimensional banded system which can be solved by QR adaptively.

\subsection{Operators on unbounded domains}\label{sec:unbounded}
We sometimes need to consider operators defined on unbounded domains, such as the Schr\"{o}dinger operator. For example, the 2D Davies' operator of complex harmonic oscillator \cite{dav1}
\begin{align*}
\mathcal{P}^{D} u = -h^2\Delta u + i(x^2 + y^2)u,
\end{align*}
which acts on $ L^2(\mathbb{R}^2)$. For such a problem, we have to ditch the Legendre polynomials and take the Hermite functions as the basis \cite{col1}. After reordering the basis functions in the tensor product of the Hermite functions, the infinite-dimensional matrix representation of $\mathcal{P}^{D}$ is banded with a growing bandwidth. This system is again solved by adaptive QR. \cref{fig:pseDavies2D} is an illustration of the pseudospectra of Davies' operator for $h = 0.8$, which much resembles the pseudospectra of the one-dimensional version given in \cite{dav1} and \cite[\S 43]{tre5}.

\section{Conclusion}\label{sec:closing}
We have proposed a continuous framework for computing the pseudospectra of linear operators, provided that the linear operator has a compact or compact-plus-scalar resolvent. The new method is free of spectral pollution and invisibility and features adaptivity in both the stopping criterion for the operator Lanczos iteration and the DoF for resolved solutions to the inverse resolvent equations. The adoption of the well-conditioned spectral methods prevents any deterioration in the conditioning of the problem. If a small enough tolerance is used in the exit condition of the Lanczos process, the computed resolvent norm has a near-optimal accuracy.

We are pleased to see what was once anticipated by Nick Trefethen a quarter century ago has now come to pass.

\section*{Acknowledgments}

We are grateful to Matthew Colbrook (Cambridge) and Anthony P. Austin (NPS) for sharing with us their insights and valuable feedback on a draft of this paper which led us to significantly improve our work. We have also benefited from discussion with our team members Lu Cheng, Ouyuan Qin, and Kaining Deng. Finally, we would like to thank two constructive anonymous referees.

\bibliographystyle{siamplain}
\bibliography{references}
\end{document}